\documentclass[12pt,reqno]{amsart}
\usepackage{amsmath}
\usepackage{amssymb}
\usepackage{amsthm}
\usepackage{latexsym}
\usepackage{amssymb}

\newtheorem{theorem}{Theorem}[section]
\newtheorem{lemma}{Lemma}[section]
\newtheorem{corollary}{Corollary}[section]
\newtheorem{remark}{Remark}[section]
\newtheorem{definition}{Definition}[section]
\newtheorem{proposition}{Proposition}[section]
\newtheorem{example}{Example}[section]
\newtheorem{assumption}{Assumption}[section]
\numberwithin{equation}{section}
\newcommand{\bth}{\begin{theorem}}
	\newcommand{\ethe}{\end{theorem}}

\newcommand{\bre}{\begin{remark}}
	\newcommand{\ere}{\end{remark}}

\newcommand{\ble}{\begin{lemma}}
	\newcommand{\ele}{\end{lemma}}
\newcommand{\bde}{\begin{definition}}
	\newcommand{\ede}{\end{definition}}

\newcommand{\bco}{\begin{corollary}}
	\newcommand{\eco}{\end{corollary}}

\newcommand{\bpr}{\begin{proposition}}
	\newcommand{\epr}{\end{proposition}}

\newcommand{\bexer}{\begin{exercise}}
	\newcommand{\eexer}{\end{exercise}}
\newcommand{\breh}{\begin{hint}}
	\newcommand{\ereh}{\end{hint}}

\newcommand{\halmos}{\hfill \qed}

\newcommand{\bexam}{\begin{example}}
	\newcommand{\eexam}{\end{example}}

\newcommand{\pr} {{\bf Proof.}}

\newcommand{\bfi}{\begin{fig}}
	\newcommand{\efi}{\end{fig}}

\newcommand{\beao}{\begin{eqnarray*}}
	\newcommand{\eeao}{\end{eqnarray*}\noindent}

\newcommand{\beam}{\begin{eqnarray}}
	\newcommand{\eeam}{\end{eqnarray}\noindent}
\newcommand{\E}{\mathbf{E}}
\newcommand{\PP}{\mathbf{P}}

\newcommand{\xto}{x\to\infty}

\newcommand{\bH}{\overline{H}}
\newcommand{\bF}{\overline{F}}

\newcommand{\bG}{\overline{G}}

\newcommand{\bbr}{{\mathbb R}}

\newcommand{\bbn}{{\mathbb N}}

\newcommand{\vep}{\varepsilon}

\usepackage{soul}
\setstcolor{red}

\usepackage{color}

\begin{document}
	\title{Tail behavior of randomly weighted sums with interdependent summands}

	\author[ D.G. Konstantinides, R. Leipus, C.D. Passalidis, J. \v{S}iaulys]{ Dimitrios G. Konstantinides, Remigijus Leipus, \\ Charalampos D. Passalidis, Jonas \v{S}iaulys}

	\address{Dept.\ of Statistics and Actuarial-Financial Mathematics,
		University of the Aegean,
		Karlovassi, Samos, Greece}
	\email{konstant@aegean.gr}
	
	\address{Inst.\ of Applied Mathematics,
		Vilnius University,
		Vilnius, Lithuania}
	\email{remigijus.leipus@mif.vu.lt}
	
	\address{Dept.\ of Statistics and Actuarial-Financial Mathematics,
		University of the Aegean,
		Karlovassi, Samos, Greece}
	\email{sasm23002@sas.aegean.gr}
	
	\address{Inst.\ of Mathematics,
		Vilnius University,
		Vilnius, Lithuania}
	\email{jonas.siaulys@mif.vu.lt}

	\date{{\small \today}}
	
	\begin{abstract} 
		We reconsider a classical, well-studied problem from applied probability. This is the max-sum equivalence of randomly weighted sums, and the originality is because we manage to include interdependence among the primary random variables, as well as among primary random variables and random weights, as a generalization of previously published results. As a consequence we provide the finite-time ruin probability, in a discrete-time risk model.  Furthermore, we established asymptotic bounds for the generalized moments of randomly weighted sums in the case of dominatedly varying primary random variables under the same dependence conditions. Finally, we give some results for randomly weighted and stopped sums under similar dependence conditions, with the restriction that the random weights are identically distributed, and the same holds for the primary random variables. Additionally, under these assumptions, we find asymptotic expressions for the random time ruin probability, in a discete-time risk model.
	\end{abstract}
	
	\maketitle
	\vspace{3mm}
	

	\section{Introduction} \label{sec.KLP.1}
	
	Let $X_1,\,\ldots,\,X_n$, $n\in \bbn$, be real-valued random variables (r.v.s), with distributions  
	$F_1,\,\ldots,\,F_n$, respectively, and let $\Theta_1,\,\ldots,\,\Theta_n$ be non-negative non-degenerate at zero r.v.s, with distributions $G_1,\,\ldots,\,G_n$, respectively. The r.v.s $X_i$, for $i=1,\,\ldots,\,n$, represent the primary random variables and the r.v.s $\Theta_i$, for $i=1,\,\ldots,\,n$, depict the random weights. In this paper, we consider the weighted sum 
	\beam \label{eq.KLP.1}
	S_n^{\,\Theta}&:=&\sum_{i=1}^n \Theta_i X_i\,,
	\eeam
	and study the asymptotic behavior of the tail probability of this sum and the corresponding maximum of such sums and the maximal summands, defined below:
	\beao
	\bigvee_{i=1}^n S_i^{\,\Theta}&:=&\max_{1\leq k \leq n} \sum_{i=1}^k \Theta_i X_i\,,\qquad \bigvee_{i=1}^n \Theta_i X_i\ \ :=\ \ \max_{1\leq i \leq n} \Theta_i X_i.
	\eeao
	Therefore, a crucial step is the proof of relation
	\beam \label{eq.KLP.3}
	\PP\left[ S_n^{\,\Theta}>x \right] &\sim& \sum_{i=1}^ n \PP[\Theta_i \,X_i >x]\,,
	\eeam
	as $\xto$, or furthermore, the proof of relations
	\beam \label{eq.KLP.4}
	\PP\left[ S_n^{\,\Theta}>x \right]  &\sim&\PP\left[ \bigvee_{i=1}^n S_i^{\,\Theta}>x \right] \ \ \sim\ \ \PP\left[ \bigvee_{i=1}^n \Theta_i\,X_i>x \right] \ \ \sim\ \  \sum_{i=1}^ n \PP[\Theta_i \,X_i >x],
	\eeam
	as $\xto$.
	
	The asymptotic relations \eqref{eq.KLP.4} were studied in the frame of heavy-tailed distributions (usually, subexponential distributions, see definition in Section \ref{sec.KLP.2} below) with applications mainly in risk theory, where the $X_i$ represent claim sizes. At the same time, the $\Theta_i$ plays the role of the stochastic discount factor, while in the case of degenerate to constant nonzero value, the $\Theta_i$ depicts the discount factor in a model with a constant force of interest. In most cases, relations \eqref{eq.KLP.3} and \eqref{eq.KLP.4} were considered under the condition that $\{ \Theta_1,\,\ldots,\,\Theta_n \}$ are arbitrarily dependent r.v.s, independent of  $\{ X_1,\,\ldots,\,X_n \}$, with the last also mutually independent, see for example \cite{tang:tsitsiashvili:2003b} and \cite{tang:yuan:2014}.
	
	The first attempts to study this problem under dependence structures, that contain independence as a particular case, were oriented to the non-weighted sum 
	\beao
	S_n := \sum_{i=1}^n X_i,
	\eeao 
	as e.g.\ papers \cite{geluk:tang:2009} and \cite{ko:tang:2018}.
	
	In the sequence, several papers were published considering the asymptotic behavior of the tail of $S_n^{\,\Theta}$ in the following two cases:
		\begin{enumerate}
			\item r.v.s 
			$X_1,\,\ldots,\,X_n$  are somehow dependent, while r.v.s $\Theta_1,\,\ldots,\,\Theta_n$ are arbitrarily dependent and collections $\{X_1,\,\ldots,\,X_n\}$,  $\{\Theta_1,\,\ldots,\,\Theta_n\}$ are mutually independent,
			\item
			primary r.v.\ $X_i$ and weight $\Theta_i$ are somehow dependent for any index $i$, and the pairs $(X_i,\,\Theta_i)$, $(X_j,\,\Theta_j)$ are independent for any $i\neq j$. 
	\end{enumerate}
	For (1) see, e.g., \cite{cheng:2014} and \cite{li:2013}, while for (2) see, e.g., \cite{yang:leipus:siaulys:2012} and \cite{yang:wang:leipus:siaulys:2013}. Taking into account that in practical applications, one can find the cases where both (1) and (2) appear simultaneously, it is of interest to generalize some already published works by combining (1) and (2) through a new dependence structure, inspired by \cite{asimit:jones:2008} and next by the generalization to multidimensional setup in \cite{jiang:wang:chen:xu:2015} and \cite{li:2016}. In the literature, (1) and (2) were combined only in the case of regularly varying distributions under the framework of multivariate regular variation and some multivariate extensions of Breiman's lemma, see for example \cite{fougeres:mercadier:2012}. We refer to the combination of (1) and (2) cases as {\em interdependence}. We would like to remind that, according to our knowledge, the only paper, where this kind of 'interdependence' is contained with distribution classes larger than regular variation, is the paper \cite{chen:cheng:2024}, where the relation \eqref{eq.KLP.3} was examined in the case of $n=2$, with some dependence structures inspired by \cite{asimit:jones:2008}, similar to our structures. 
	
	The rest of the paper is organized as follows. In Section \ref{sec.KLP.2}, after some basic definitions, we provide some heavy-tailed distribution classes and some dependence structures that play an essential role in our work. In Section \ref{sec.KP.2}, we introduce a new dependence structure and present the main results (related with \eqref{eq.KLP.4}) with their proofs after some preliminary lemmas. In Section \ref{sec.KP.3}, we proceed to the partial extension of the results in case of generalized moments of randomly weighted sums. In Section \ref{sec.KP.5}, we study the max-sum equivalence of randomly weighted stopped sums under some more restrictive conditions. It is worth mentioning that in sections \ref{sec.KP.2} and  \ref{sec.KP.5}, the asymptotic estimates for the ruin probability of finite and random (finite) time in the discrete-time risk model are given. Meanwhile, the asymptotic estimates of generalized moments, from Section \ref{sec.KP.3}, permit direct applications to risk measures, such as expected shortfall and marginal expected shortfall.
	
	\section{Preliminaries}\label{sec.KLP.2}
	
	\subsection{Notations}
	
	Throughout the paper, we consider the real-valued r.v.s $X_1,\dots,X_n$ with the corresponding distributions $F_1,\dots,F_n$ as primary r.v.s., while 
	$\Theta_1,\dots\Theta_n$, with corresponding distributions $G_1,\dots,G_n$, are random weights, which are non-negative non-degenerate at zero r.v.s. 
	
	Further, we denote by $\bF_i:=1-F_i$ the distribution tail, hence $\bF_i(x)=\PP[X_i>x]$, and assume that $\bF_i(x)>0$ for any $x \geq 0$ and   $i=1,\,\ldots,\,n$, although this is not necessarily true for distributions $\bG_i(x)$. For two positive functions $f(x)$ and $g(x)$, the asymptotic relation $f(x)=o[g(x)]$, as $\xto$, denotes the convergence
	\beao
	\lim_{\xto}\dfrac{f(x)}{g(x)} = 0\,,
	\eeao
	while the asymptotic relation $f(x)=O[g(x)]$, as $\xto$, means
	\beao
	\limsup_{\xto} \dfrac{f(x)}{g(x)} < \infty,\eeao
	and the asymptotic relation $f(x)\asymp g(x)$, as $\xto$ means that hold both $f(x)=O[g(x)]$ and 
	$g(x)=O[f(x)]$. Further we write $f(x)\gtrsim g(x)$ or $f(x)\lesssim g(x)$ if it holds that
	\beao
	\liminf_{\xto}\frac{f(x)}{g(x)}\geq 1\,, \;\;\text{or } \;\; \limsup_{\xto}\frac{f(x)}{g(x)}\leq 1\,,
	\eeao
	respectively. For a real number $x$, we denote $x^+:=\max\{x,0\}$
	and $x^-:=-\min\{x,0\}$. For two real numbers $x$ and $y$ we denote $x\vee y := \max\{x,\,y\}$ and $x \wedge y:=\min\{x,\,y\}$. For some event $A$, the symbol ${\bf 1}_A$ depicts the indicator function of $A$. 
	
	\subsection{Heavy tailed distributions}
	
	Next, we consider some heavy-tailed distribution classes. We say that a distribution $F$ is heavy-tailed, symbolically $F \in \mathcal{H}$, if
	\beao
	\int_{-\infty}^{\infty} e^{\vep\,y}F(dy)=\infty,
	\eeao
	for any $\vep >0$. We say that a distribution $F$ has a long tail, symbolically $F \in \mathcal{L}$, if  for any (or equivalently, for some) $a> 0$ it holds that
	\beao
	\lim_{\xto} \dfrac{\bF(x-a)}{\bF(x)}  = 1.
	\eeao
	
	Next distribution class represents a popular subclass of heavy-tailed distributions, with several applications in risk theory, financial mathematics, branching processes and queuing theory, see for example \cite{asmussen:2003}, \cite{asmussen:steffensen:2020}, \cite{foss:korshunov:zachary:2013} and \cite{Konstantinides:2018}. We say that a distribution $F$ with support $[0,\,\infty)$ belongs to the subexponential class of distributions, symbolically $F \in \mathcal{S}$ if for any (or equivalently, for some) integer $n\geq 2$, it holds that
	\beao
	\lim_{\xto} \dfrac{\overline{F^{*n}}(x)}{\bF(x)}  = n\,,
	\eeao
	where $F^{*n}$ depicts the $n$-th order convolution power of $F$. Furthermore, for some $F$ with support the whole real axis $\bbr$, we say that $F\in \mathcal{S}$ if $F^+ \in \mathcal{S}$, where $F^+(x):=F(x)\,{\bf 1}_{[0,\infty)}(x)$, see for example \cite[Ch.\ 3.2]{foss:korshunov:zachary:2013}.
	
	Another known class of heavy-tailed distributions is the class of dominatedly varying distributions, introduced in \cite{feller:1969}. We say that a distribution $F$ has dominated variation, symbolically $F \in \mathcal{D}$, if for any (or equivalently, for some) $b \in (0,\,1)$ it holds
	\beao
	\limsup_{\xto} \dfrac{\bF(b\,x)}{\bF(x)}  < \infty\,.
	\eeao
	It is well-known that $\mathcal{D} \not\subseteq \mathcal{S}$, $\mathcal{S} \not\subseteq \mathcal{D}$ and $\mathcal{D}\cap \mathcal{S}=\mathcal{D}\cap \mathcal{L}\neq \varnothing$, see for example \cite{goldie:1978}.
	
	Another class of heavy-tailed distributions is the class of consistently varying distributions. We say that a distribution $F$ is consistently varying, symbolically $F \in \mathcal{C}$, if it holds
	\beao
	\lim_{v\uparrow 1} \limsup_{\xto} \dfrac{\bF(v\,x)}{\bF(x)}  =1.
	\eeao
	
	One even smaller class, with high popularity in applied probability, is the class $\mathcal{R}$ of the regularly varying distributions. We say that the distribution $F$ is regularly varying, with some index $0<\alpha <\infty$, symbolically $F\in \mathcal{R}_{-\alpha}$, if it holds
	\beao
	\lim_{\xto} \dfrac{\bF(b\,x)}{\bF(x)} = b^{-\alpha}\,,
	\eeao
	for any $b>0$.
	It is well-known that $$\mathcal{R}:=\bigcup_{0<\alpha < \infty} \mathcal{R}_{-\alpha}\subsetneq  \mathcal{C} \subsetneq \mathcal{D}\cap \mathcal{L} \subsetneq  \mathcal{S} \subsetneq  \mathcal{L} \subsetneq  \mathcal{H},$$ see for example \cite{foss:korshunov:zachary:2013} and \cite{leipus:siaulys:konstantinides:2023}.
	
	Now, we present the Matuszewska indexes and the $L$-indexes. For any distribution $F$, with condition $\bF(x)>0$ for all $x>0$ and for any $v>1$, we define 
	\beao
	\bF_{*}(v):=\liminf_{x\rightarrow \infty }\dfrac{\bF(v\,x)}{\bF(x)}\,, \qquad \bF^{*}(v):=\limsup_{x\rightarrow \infty }\dfrac{\bF(v\,x)}{\bF(x)}.
	\eeao
	Next, we bring up the upper and lower Matuszewska indexes, introduced in \cite{matuszewska:1964}:
	\beao
	\alpha_{F} :=\inf\left\{- \dfrac{\ln \bF_{*}(v)}{\ln v}\;:\;v>1 \right\}= -\lim_{y \rightarrow \infty }\dfrac{\ln \bF_{*}(y )}{\ln y }\,, \\[2mm]
	\beta_{F} :=\sup \left\{- \dfrac{\ln \bF^{*}(v)}{\ln v}\;:\;v>1 \right\}= -\lim_{y \rightarrow \infty }\dfrac{\ln \bF^{*}(y )}{\ln y }\,.
	\eeao
	It holds that $0\leq \beta_{F} \leq \alpha_{F}\le\infty$. We have the following equivalence: $F\in \mathcal{D}$ if and only if $\alpha_{F} < \infty$. 
	
	Two other popular indexes are the $L$-indexes, defined for distributions $F$, with condition $\bF(x)>0,\, x>0$. In this case, we define
	\beam \label{eq.KLP.5}
	L_F:=\lim_{v \downarrow 1} \bF_{*}(v)\,,\quad L^F:=\lim_{v \downarrow 1} \bF^{*}(v)\,,
	\eeam
	whence from relation \eqref{eq.KLP.5} we find out $0\leq L_F \leq L^F \leq 1$. Furthermore, the inclusion $F \in \mathcal{D}$ is equivalent to the inequality $L_F>0$, and the inclusion $F \in \mathcal{C}$ is equivalent to the equality $L_F=1$, see for example, \cite[Sec.\ 2.4]{leipus:siaulys:konstantinides:2023}. Several applications of index $L_F$ can be found, e.g., in \cite{cheng:2014} and \cite{wang:wang:cheng:2006}.
	
	\subsection{Dependence modeling}
	
	Now we introduce two basic dependence structures, used to describe the dependence among the primary r.v.s $X_1,\,\ldots,\,X_n$, containing the independence as a special case.

        $\bullet$ \textit{The real r.v.s $X_1,\,\ldots,\,X_n$, with distributions $F_1,\,\ldots,\,F_n$ respectively, are said to be pairwise quasi-asymptotically independent, symbolically $pQAI$, if for any pair $i,\,j$ with $i\neq j$ it holds that
	\beao
	\lim_{\xto} \dfrac{\PP[X_i^+>x\,,\;X_j^+ >x]}{\bF_i(x) + \bF_j(x)}=\lim_{\xto} \dfrac{\PP[X_i^->x\,,\;X_j >x]}{\bF_i(x) + \bF_j(x)}=0\,,
	\eeao
	or equivalently $\lim_{\xto} \PP\left[|X_i|\wedge X_j>x\;\big|\;X_i\vee X_j >x\right]=0$. }
   
$\bullet$ \textit{	We say that $X_1,\dots,X_n$ are pairwise tail asymptotically independent, symbolically $pTAI$, if for any pair $i,\,j\in\{1,2,\ldots,n\}$ with $i\neq j$ it holds that
	\beao
	\lim_{x_i \wedge x_j \to \infty} \PP\left[|X_i|>x_i\;\big|\; X_j >x_j\right]=0\,.
	\eeao}
    
	These dependence structures were introduced in \cite{chen:yuen:2009} and \cite{geluk:tang:2009} respectively. When there are only two random variables, we simply write $QAI$ and $TAI$. These dependence structures were studied in several papers, for example, \cite{chen:liu:2022}, \cite{cheng:2014}, \cite{leipus:paukstys:siaulys:2021}, \cite{li:2013} and \cite{yang:wang:konstantinides:2014}.
	
	The following known result contains two parts. The first one was initially established for non-weighted random sums in Th.\ 3.1 of \cite{geluk:tang:2009} and later for weighted sums with random weights bounded from above, in Th.\ 1 of \cite{wang:2011}. The second part can be found in \cite[Th.\ 3.1]{chen:yuen:2009}.
	
	\bpr \label{pr.KLP.1}
	Let $Z_1,\,\ldots,\,Z_n$ be real-valued r.v.s with distributions  $H_1,\,\ldots,\,H_n$ respectively.
	\begin{enumerate}
		\item
		If $H_i \in \mathcal{D} \cap \mathcal{L}$ for  $i=1,\,\ldots,\,n$ and $Z_1,\,\ldots,\,Z_n$ are $pTAI$, then
		\beao
		\PP\left[ \sum_{i=1}^{n} Z_i>x \right] \sim\PP\left[ \bigvee_{i=1}^n Z_i>x \right] \sim \PP\left[ \bigvee_{k=1}^n \sum_{i=1}^k Z_i>x \right] \sim \sum_{i=1}^ n \PP[Z_i>x]\,,
		\eeao
		as $\xto$.
		\item
		If $H_i \in \mathcal{C}$ for  $i=1,\,\ldots,\,n$ and $Z_1,\,\ldots,\,Z_n$ are $pQAI$, then
		\beao
		\PP\left[ \sum_{i=1}^{n} Z_i>x \right] \sim \sum_{i=1}^ n \PP[Z_i >x]\,,
		\eeao
		as $\xto$.
	\end{enumerate}
	\epr
	
	In this paper, we generalize Proposition \ref{pr.KLP.1}, using random weights that have a specific dependence structure with the primary random variables and under a new dependence structure. The main results of this article, both directly related to the asymptotics of a fixed number of randomly weighted sums and to the application of the resulting asymptotic formulas, are formulated in theorems \ref{th.KLP.1}, \ref{th.KLP.2}, \ref{th.KLP.3} and \ref{th.KLP.5.1}. Several lemmas are used in the proofs, which, in our opinion, are also useful in their own way.
	
	\section{Results on randomly weighted sums} \label{sec.KP.2}
	
	This section examines the asymptotic behavior of randomly weighted sums under $pTAI$ and $pQAI$ dependence structures for the primary r.v.s and simultaneously under some dependence structure among the random weights and primary r.v.s. 
	We present now the first assumption, which is needed further in the paper.
	\begin{assumption} \label{ass.KLP.2}
	Let $X_1,\ldots, X_n$ be the primary r.v.s, and let  $\Theta_1,\ldots,\Theta_n$ be the random weights.
Suppose that there exist  functions     $b_i:\;[0,\,\infty) \longrightarrow (0,\,\infty),\, i\in\{1,\ldots,n\}$,\, such that $b_i(x) \rightarrow \infty,\ \ \lim_{\xto} b_i(x)/x = 0$,	and 
		\beam \label{eq.KLP.11}
		\bG_i[b_i(x_i)]:=\PP\left[\Theta_i>b_i( x_i)\right]=o\left(\PP[\Theta_i\,X_i >  x_i\,,\;\Theta_j\,X_j > x_j]\right)\,,
		\eeam 
		as $x_i\wedge x_j \rightarrow \infty$, for $1\leq i \neq j \leq n$.
	\end{assumption} 
	
	\bre \label{rem.KLP.2} 
	Let $i,j$ be two  indices $1\le i\neq j\le n$. If r.v.s $X_i$ and $X_j$ have supports with infinite right endpoints, then it is easy to see that \eqref{eq.KLP.11} holds in the case where $G_i$ or $G_j$  have upper bounded support. In addition, we observe that relation \eqref{eq.KLP.11} implies that 
	\beao
	\dfrac{\bG_i[b_i(x_i)]}{\PP[\Theta_i\,X_i >  x_i]}\leq \dfrac{\bG_i[b_i(x_i)]}{\PP[\Theta_i\,X_i >  x_i\,,\;\Theta_j\,X_j > x_j]} \rightarrow 0\,,
	\eeao
	as $x_i\wedge x_j \rightarrow \infty$. Hence, Assumption \ref{ass.KLP.2} implies that  
	\beao
	\bG_i[b_i(x_i)]=o(\PP[\Theta_i\,X_i >  x_i])\,,
	\eeao
	as $x_i\rightarrow \infty$. In combination with properties of functions $b_i$ and \cite[Lemma  3.2]{tang:2006} we find
	\beam \label{eq.KLP.13}
	\bG_i(c_i\,x_i)=o\left(\PP[\Theta_i\,X_i >  x_i]\right)\,,
	\eeam
	as $x_i \rightarrow \infty$, respectively for any constants $c_i > 0$. Note that relation \eqref{eq.KLP.13} can be useful in the study of closure properties with respect to convolution product; see, for example, \cite{konstantinides:leipus:siaulys:2022} and \cite{tang:2006}.
	\ere
Let us now render a new dependence structure that plays a crucial role in modeling interdependence.
\begin{assumption} \label{ass.KLP.3*}
Let $X_1,\ldots, X_n$ be the primary r.v.s, and let $\Theta_1,\ldots,\Theta_n$ be the random weights. We assume that there exists some measurable functions ${g_{i j}:[0,\,\infty)^2 \rightarrow (0,\,\infty)}$, $1\le i\neq j\le n$, such that
		\beam \label{eq.KP.13*}		\PP[|X_i|>x_i\,,\;X_j>x_j\;\big|\;\Theta_i=\theta_i\,,\;\Theta_j=\theta_j]\sim g_{i j}(\theta_i,\,\theta_j)\,\PP[|X_i|>x_i\,,\;X_j>x_j]\,,
		\eeam
		as $x_i\wedge x_j\rightarrow \infty$, with $1\leq i \neq j \leq n$, uniformly for any $(\theta_i,\,\theta_j)\in s_{\Theta_i}\times s_{\Theta_j}$, where $s_{\Theta_i}$ represents the support of $\Theta_i$, for  any $i=1,\,\ldots,\,n$, i.e. 
\beao
		s_{\Theta_i}:=\big\{\theta_i \in \mathbb{R}\;:\; \PP(\Theta_i\in[\theta_i-\delta,\theta_i+\delta])>0\ {\rm for\ any\ } \delta>0\big\}\,.
\eeao
	\end{assumption}	
	In this assumption, the uniformity is understood in the sense that
	\beao
	\lim_{x_i\wedge x_j) \rightarrow \infty} \sup_{(\theta_i,\,\theta_j) \in (s_{\Theta_i},\,s_{\Theta_j})}\left|\dfrac{\PP\left[ |X_i|>x_i\,,\;X_j>x_j\;\big|\;\Theta_i=\theta_i\,,\;\Theta_j=\theta_j\right]}{g_{i j}(\theta_i,\,\theta_j)\, \PP[ |X_i|>x_i\,,\;X_j>x_j]} - 1 \right|=0,
	\eeao
	where $1\le i\ne j\le 2$, and
    the conditional probability in \eqref{eq.KP.13*}	is defined by the equality
    \smallskip
    \beao
\lim_{\delta\downarrow 0}\,\frac{\PP\left[|X_i|>x_i\,,\;X_j>x_j,\,\Theta_i\in[\theta_i-\delta,\theta_i+\delta], \Theta_j\in[\theta_j-\delta,\theta_j+\delta]\right]}{\PP\left[\Theta_i\in[\theta_i-\delta,\theta_i+\delta], \Theta_j\in[\theta_j-\delta,\theta_j+\delta]\right]}\,.
    \eeao
	
	\bre \label{rem.KLP.1}
We observe that Assumption \ref{ass.KLP.3*} implies that 
\beam \label{eq.KLP.13b}
\E[g_{i j}(\Theta_i,\,\Theta_j)]=1\,,
\eeam
for any pair $1\le i\ne j\le n$.
\ere

\begin{proof}[\bf Proof] Let us denote temporally ${\bf \Theta}=(\Theta_i,\Theta_j)$. By Assumption \ref{ass.KLP.3*}, via integration on the left side of \eqref{eq.KP.13*}, we find 
	\beao
	&&\int_{[0,\infty)} \int_{[0,\infty)} {\PP[| X_{i}|>x_{i},\;X_j>x_j\;\big|\;{\bf \Theta}= (\theta_i,\,\theta_j)]}\PP[{\bf \Theta}\in(d\theta_i,\,d\theta_j)]\\[2mm]
	&&\qquad \qquad =\PP[| X_{i}|>x_{i},\;X_j>x_j]\,,
	\eeao
	and through integration on the right side of \eqref{eq.KP.13*} we obtain that
	\beao
	&&\PP[| X_{i}|>x_{i},\;X_j>x_j]\int_{[0,\infty)}  \int_{[0,\infty)} {g_{ij}(\theta_i,\,\theta_j)}\PP[{\bf \Theta}\in(d\theta_i,\,d\theta_j)]\\[2mm]
	&&\qquad \qquad =\PP[| X_{i}|>x_{i},\;X_j>x_j]\,\E\left[g_{ij}({\bf \Theta})\right]\,.
	\eeao
	Since the convergence is uniform in $(\theta_1,\,\theta_2)$, we obtain the relation
	\beao
	\PP[| X_{i}|>x_{i},\;X_j>x_j]\sim\PP[ X_{i}|>x_{i},\;X_j>x_j]\,\E\left[g_{ij}({\bf \Theta})\right]\,,
	\eeao 
	as $x_i\wedge x_j \rightarrow  \infty$, which implies $\E\left[g_{ij}({\bf \Theta})\right]=1$.
    \end{proof}

    \bre\label{aa}
We observe that in the case of non-negative primary r.v.s $X_1,\ldots, X_n$, the dependence structure described in Assumption \ref{ass.KLP.3*} is a slight generalization of Condition 2 in \cite{jiang:wang:chen:xu:2015}. We also observe that the dependence structure in Assumption \ref{ass.KLP.3*} allows for arbitrary dependence among the components of $\{X_1,\ldots, X_n\}$ and among the components of $\{\Theta_1,\ldots,\Theta_n\}$. However, it includes some form of asymptotic independence between $\{X_1,\ldots, X_n\}$ and $\{\Theta_1,\ldots,\Theta_n\}$. Further results on asymptotic idependence can be found in \cite{hazra:maulik:2012} and \cite{maulik:resnick:2004}.
	\ere 
In the following lemma, we present a result that holds for any pair $(i,\,j)$ with $1\leq i \neq j \leq n$, but for the sake of simplicity, we formulate it with indexes $1,\,2$ only. 
	
	\ble \label{lem.KLP.1} 
Let $X_1, X_2$ be primary r.v.s and let $\Theta_1, \Theta_2$ be the random weights satisfying  Assumptions \ref{ass.KLP.2} and \ref{ass.KLP.3*} for $n=2$.  Then
	\beam \label{eq.KLP.14}
	&&\PP[\Theta_1|X_1| >  x_1\,,\;\Theta_2\,X_2 > x_2]\\[2mm] \notag
	&&\qquad \sim \int_{[0,\infty)} \int_{[0,\infty)} g_{1 2}(\theta_1,\,\theta_2)\,\PP\left[|X_1| >  \dfrac{x_1}{\theta_1}\,,\;X_2 > \dfrac{x_2}{\theta_2}\right]\,\PP[\Theta_1 \in d\theta_1\,,\;\Theta_2 \in d\theta_2]\,,
	\eeam
	as $x_1\wedge x_2 \rightarrow \infty$.
	\ele
	
	\pr~
	Using the functions $b_1$ and $b_2$ from Assumption \ref{ass.KLP.2}, we obtain
	\beao
	&&\PP[\Theta_1\,|X_1| >  x_1\,,\;\Theta_2\,X_2 > x_2]\\[2mm] 
	&&= \int_{[0,\infty)} \int_{[0,\infty)} \PP\left[|X_1| >  \dfrac{x_1}{\theta_1}\,,\;X_2 > \dfrac{x_2}{\theta_2}\;\Big|\;\Theta_1=\theta_1\,,\;\Theta_2=\theta_2\right]\,\PP[\Theta_1 \in d\theta_1\,,\;\Theta_2 \in d\theta_2]\\[2mm] 
	&&= \left(\int_{[0,b_1(x_1)]}+\int_{(b_1(x_1),\infty)} \right)\,\left(\int_{[0,b_2(x_2)]}+\int_{(b_2(x_2),\infty)} \right) \\[2mm] 
	&&\qquad \PP\left[|X_1| >  \dfrac{x_1}{\theta_1}\,,\;X_2 > \dfrac{x_2}{\theta_2}\;\Big|\;\Theta_1=\theta_1\,,\;\Theta_2=\theta_2\right]\,\PP[\Theta_1 \in d\theta_1\,,\;\Theta_2 \in d\theta_2]\\[2mm] 
	&&=I_{11}(x_1,\,x_2)+I_{12}(x_1,\,x_2)+I_{21}(x_1,\,x_2)+I_{22}(x_1,\,x_2)\,.
	\eeao
	By \eqref{eq.KLP.11} we get
	\beam \label{eq.KLP.16} \notag
	I_{22}(x_1,\,x_2) &\leq & \PP\left[\Theta_1 >  b_1(x_1)\,,\;\Theta_2 > b_2(x_2)\right]\\[2mm]
	&\leq& \PP\left[\Theta_1 >  b_1(x_1)\right]=o\left(\PP\left[\Theta_1\,X_1 >  x_1\,,\;\Theta_2\,X_2 > x_2\right]\right)\\ \notag
	&=&o\left(\PP\left[\Theta_1\,|X_1| >  x_1\,,\;\Theta_2\,X_2 > x_2\right]\right)\,,
	\eeam
	as $x_1\wedge x_2 \rightarrow \infty$. Similarly,  
	\beam \label{eq.KLP.17} \notag
	I_{21}(x_1,\,x_2)&=&\int_{(b_1(x_1),\infty)} \int_{[0,b_2(x_2)]}  \PP\left[|X_1| >  \dfrac{x_1}{\theta_1}\,,\;X_2 > \dfrac{x_2}{\theta_2}\;\Big|\;\Theta_1=\theta_1\,,\;\Theta_2=\theta_2\right]\,\\[2mm] \notag
	&&\qquad \times \PP[\Theta_1 \in d\theta_1\,,\;\Theta_2 \in d\theta_2] \leq \PP\left[\Theta_1 >  b_1(x_1)\,,\;\Theta_2 \leq b_2(x_2)\right]\\[2mm]
	&\leq& \PP\left[\Theta_1 >  b_1(x_1)\right]=o\left(\PP\left[\Theta_1\,X_1 > x_1\,,\;\Theta_2\,X_2 >x_2\right]\right)\\[2mm] \notag
	&=&o\left(\PP\left[\Theta_1\,|X_1| >  x_1\,,\;\Theta_2\,X_2 > x_2\right]\right)\,,
	\eeam 
	as $x_1\wedge x_2 \rightarrow \infty$, where the second last step follows by \eqref{eq.KLP.11}. By symmetry, we have
	\beam \label{eq.KLP.18} 
	I_{12}(x_1,\,x_2)=o\left(\PP\left[\Theta_1\,X_1 > x_1\,,\;\Theta_2\,X_2 >x_2\right]\right)=o\left(\PP\left[\Theta_1\,|X_1| >  x_1\,,\;\Theta_2\,X_2 > x_2\right]\right)\,,
	\eeam 
	as $x_1\wedge x_2 \rightarrow \infty$. Finally, from Assumption \ref{ass.KLP.2} and Assumption \ref{ass.KLP.3*}  we derive the following asymptotic relation
	\beam \label{eq.KLP.18*}
&&I_{11}(x_1,\,x_2)=\int_{[0,b_1(x_1)]} \int_{[0,b_2(x_2)]}  \PP\left[|X_1| >  \dfrac{x_1}{\theta_1}\,,\;X_2 > \dfrac{x_2}{\theta_2}\;\Big|\;\Theta_1=\theta_1\,,\;\Theta_2=\theta_2\right]\\[2mm] \notag
	&&\qquad \qquad \times\, \PP[\Theta_1 \in d\theta_1\,,\;\Theta_2 \in d\theta_2] \\[2mm] \notag
&&\sim \int_{[0,b_1(x_1)]} \int_{[0,b_2(x_2)]} g_{1 2}(\theta_1,\,\theta_2)\PP\left[|X_1| >  \dfrac{x_1}{\theta_1}\,,\;X_2 > \dfrac{x_2}{\theta_2}\right] \PP[\Theta_1 \in d\theta_1\,,\;\Theta_2 \in d\theta_2]\,, 
	\eeam
	as $x_1\wedge x_2 \rightarrow \infty$. Indeed, for any $\epsilon$ and large $x_1\wedge x_2$ we have
	\beao
	&&\frac{\int\limits_{[0,b_1(x_1)]} \int\limits_{[0,b_2(x_2)]}  \PP\left[|X_1| >  \dfrac{x_1}{\theta_1},X_2 > \dfrac{x_2}{\theta_2}\Big|\Theta_1=\theta_1,\Theta_2=\theta_2\right]\PP[\Theta_1 \in d\theta_1,\Theta_2 \in d\theta_2]}
	{\int\limits_{[0,b_1(x_1)]} \int\limits_{[0,b_2(x_2)]} g_{1 2}(\theta_1,\theta_2)\PP\left[|X_1| > \dfrac{x_1}{\theta_1},X_2 > \dfrac{x_2}{\theta_2}\right] \PP[\Theta_1 \in d\theta_1,\Theta_2 \in d\theta_2]}
	\\
	&&\le \sup_{(\theta_1,\theta_2)\in (s_{\Theta_1},s_{\Theta_2})}\frac{\PP\left[|X_1| >  \dfrac{x_1}{\theta_1},X_2 > \dfrac{x_2}{\theta_2}\Big|\Theta_1=\theta_1,\Theta_2=\theta_2\right]}
	{g_{1 2}(\theta_1,\theta_2)\PP\left[|X_1| > \dfrac{x_1}{\theta_1},X_2 > \dfrac{x_2}{\theta_2}\right]}
	\ \le\ 1+\epsilon,
	\eeao
	and similarly 
	\beao
	&&\frac{\int\limits_{[0,b_1(x_1)]} \int\limits_{[0,b_2(x_2)]}  \PP\left[|X_1| >  \dfrac{x_1}{\theta_1},X_2 > \dfrac{x_2}{\theta_2}\Big|\Theta_1=\theta_1,\Theta_2=\theta_2\right]\PP[\Theta_1 \in d\theta_1,\Theta_2 \in d\theta_2]}
	{\int\limits_{[0,b_1(x_1)]} \int\limits_{[0,b_2(x_2)]} g_{1 2}(\theta_1,\theta_2)\PP\left[|X_1| > \dfrac{x_1}{\theta_1},X_2 > \dfrac{x_2}{\theta_2}\right] \PP[\Theta_1 \in d\theta_1,\Theta_2 \in d\theta_2]}\\[2mm]
	&&\qquad \geq \inf_{(\theta_1,\theta_2)\in (s_{\Theta_1},s_{\Theta_2})}\frac{\PP\left[|X_1| >  \dfrac{x_1}{\theta_1},X_2 > \dfrac{x_2}{\theta_2}\Big|\Theta_1=\theta_1,\Theta_2=\theta_2\right]}
	{g_{1 2}(\theta_1,\theta_2)\PP\left[|X_1| > \dfrac{x_1}{\theta_1},X_2 > \dfrac{x_2}{\theta_2}\right]}\ 
	\ge \ 1-\epsilon,
	\eeao 
	which implies the asymptotic relation \eqref{eq.KLP.18*}.
	
	The derived asymptotic relations \eqref{eq.KLP.16}--\eqref{eq.KLP.18*} 
	imply that
	\beam \label{eq.KLP.18b} \notag
	&&\PP[\Theta_1\,|X_1| >  x_1\,,\;\Theta_2\,X_2 > x_2]\\[2mm]
	&&=[1+o(1)]\,\int\limits_{[0,b_1(x_1)]} \int\limits_{[0,b_2(x_2)]}  \PP\left[|X_1| >  \dfrac{x_1}{\theta_1}\,,\;X_2 > \dfrac{x_2}{\theta_2}\;\Big|\;\Theta_1=\theta_1\,,\;\Theta_2=\theta_2\right]\\[2mm] \notag
	&&\qquad  \times\,\PP[\Theta_1 \in d\theta_1\,,\;\Theta_2 \in d\theta_2] +o\left(\PP\left[\Theta_1\,|X_1| >  x_1\,,\;\Theta_2\,X_2 > x_2\right]\right)=\\[2mm] \notag
	&&[1+o(1)] \left(\int\limits_{[0,\infty)} \int\limits_{[0,\infty)} -\int\limits_{[0,b_1(x_1)]}\int\limits_{(b_2(x_2),\infty)} -\int\limits_{(b_1(x_1),\infty)} \int\limits_{[0,b_2(x_2)]}-\int\limits_{(b_1(x_1),\infty)} \int\limits_{(b_2(x_2),\infty)}   \right) g_{1 2}(\theta_1, \theta_2) \\[2mm] \notag
	&&\times \PP\left[|X_1| >  \dfrac{x_1}{\theta_1}\,,\;X_2 > \dfrac{x_2}{\theta_2}\right]\,\PP[\Theta_1 \in d\theta_1\,,\;\Theta_2 \in d\theta_2]+\,o\left(\PP\left[\Theta_1\,|X_1| >  x_1\,,\;\Theta_2\,X_2 > x_2\right]\right)\,,
	\eeam
	as $x_1\wedge x_2 \rightarrow \infty$.
	
Since function $g_{12}(\theta_1,\theta_2)$ is bounded, the integrals 
		\begin{eqnarray*}
			&&\bigg(\int\limits_{[0,b_1(x_1)]}\int\limits_{(b_2(x_2),\infty)} +\int\limits_{(b_1(x_1),\infty)} \int\limits_{[0,b_2(x_2)]}+\int\limits_{(b_1(x_1),\infty)} \int\limits_{(b_2(x_2),\infty)}\bigg) g_{1 2}(\theta_1,\,\theta_2)\,\\[2mm] \notag
			&&\quad \quad\times \PP\left[|X_1| >  \dfrac{x_1}{\theta_1}\,,\;X_2 > \dfrac{x_2}{\theta_2}\right]\,\PP[\Theta_1 \in d\theta_1\,,\;\Theta_2 \in d\theta_2]\,,
		\end{eqnarray*}
		can be estimated by $o\left(\PP\left[\Theta_1\,|X_1| >  x_1\,,\;\Theta_2\,X_2 > x_2\right]\right)$ like in relations \eqref{eq.KLP.16}, \eqref{eq.KLP.17} and \eqref{eq.KLP.18}. Hence, the asymptotic relation 
		\eqref{eq.KLP.18b} imply \eqref{eq.KLP.14}. This finishes the proof of  Lemma \ref{lem.KLP.1}. 
	~\halmos
	
	\bre \label{rem.KLP.3}
Following the argument from  \cite[pages 3331--3332]{chen:xu:cheng:2019} we can derive that, in the case of Assumption~\ref{ass.KLP.3*}, there exists some constants $K_{ij}>0$, $1\le i\neq j\le n$, such that $g_{ij}(\theta_i,\,\theta_j) \leq K_{ij}$ for any  $(\theta_i,\,\theta_j) \in (s_{\Theta_i},\,s_{\Theta_j})$.
	\ere
	
\begin{proof}[\bf Proof]	It is sufficient to consider the case $n=2$. Let us assume at moment the contrary, namely that for any $K < \infty$, there exists a vector $(\theta_1^{\circ},\,\theta_2^{\circ}) \in (s_{\Theta_1},\,s_{\Theta_2})$, such that $g_{12}(\theta_1^{\circ},\,\theta_2^{\circ}) > K$. Then, according to  Assumption \ref{ass.KLP.3*}
	\beao
\PP\left[|X_1| > x_1,\,X_2 >x_2\;\big|\;\Theta_1=\theta_1^{\circ},\,\Theta_2=\theta_2^{\circ}\right] \sim g_{12}(\theta_1^{\circ},\,\theta_2^{\circ})\,\PP\left[|X_1| > x_1,\,X_2 >x_2\right]\,,
	\eeao
	which implies that 
	\beao
\PP\left[|X_1| > x_1,\,X_2 >x_2\;\big|\;\Theta_1=\theta_1^{\circ},\,\Theta_2=\theta_2^{\circ}\right] \gtrsim K\,\PP\left[|X_1| > x_1,\,X_2 >x_2\right]\,,
\eeao
as $x_{1}\wedge x_{2} \rightarrow  \infty$.  Therefore, the conditional probability on the left-hand side becomes greater than unity: 
\beao
\liminf_{x_{1}\wedge x_{2} \rightarrow \infty} \PP\left[|X_1| > x_1,\,X_2 >x_2\right] > \dfrac 1K\,,
\eeao 
which leads to absurd. 
\end{proof}
	
The next result shows the closure property of two dependence structures, $pQAI$ and $pTAI$, with respect to product convolution under Assumption \ref{ass.KLP.2}, for distributions $F_1,\,\ldots,\,F_n$ from the class $\mathcal{D}$ and under independence between $\{ X_1,\,\ldots,\,X_n \}$ and $\{ \Theta_1,\,\ldots,\,\Theta_n \}$. A similar result under different conditions on $\Theta_1,\,\ldots,\,\Theta_n$ can be found in \cite[Th.\ 2.2]{li:2013}. Those conditions based on the finite moments of $\Theta_1,\,\ldots,\,\Theta_n$, of some order greater than $\max_{1\leq i \leq n} \alpha_{F_i} $, the maximum of upper Matuszewska indexes of $F_i$.
	
\ble \label{lem.KLP.2}
If $X_1,\,\ldots,\,X_n$ are $pTAI$ (or $pQAI$) real-valued r.v.s, with distributions  $F_1,\,\ldots,\,F_n \in \mathcal{D}$ respectively, and $\Theta_1,\,\ldots,\,\Theta_n$ are non-negative, non-degenerate at zero r.v.s, independent of $X_1,\,\ldots,\,X_n$. Under Assumption  \ref{ass.KLP.2} 
the products $\Theta_1 X_1,\,\ldots,\,\Theta_n X_n$ are also $pTAI$ (or $pQAI$, respectively) random variables.
	\ele 
	
	\pr~
	Let us begin with the $pTAI$ case. For positive $x_i>0$ and $x_j>0$, with $1\leq i \neq j \leq n$ we write
	\begin{align} \label{eq.KLP.21} \notag
		\PP\left[|\Theta_i\,X_i| >  x_i\,,\;\Theta_j\,X_j > x_j\right] &=J_{11}(x_i,x_j) + J_{12}(x_i,x_j) \nonumber\\&\quad+ J_{21}(x_i,x_j) + J_{22}(x_i,x_j),
	\end{align} 
	where
	\beao
	J_{11}(x_i,x_j) &:=& \int_{[0,b_i(x_i)]} \int_{[0,b_j(x_j)]} 
	\PP\left[|X_i|>\dfrac {x_i}{\theta_{i}},\,X_j >\dfrac {x_j}{\theta_{j}}\right] \PP\left[\Theta_{i}\in d\theta_{i}, \Theta_{j}\in d\theta_{j}\right],\\
	J_{12}(x_i,x_j) &:=& \int_{[0,b_i(x_i)]} \int_{(b_j(x_j),\infty)} 
	\PP\left[|X_i|>\dfrac {x_i}{\theta_{i}},\,X_j >\dfrac {x_j}{\theta_{j}}\right] \PP\left[\Theta_{i}\in d\theta_{i}, \Theta_{j}\in d\theta_{j}\right],\\
	J_{21}(x_i,x_j) &:=& \int_{(b_i(x_i),\infty)} \int_{[0,b_j(x_j)]}
	\PP\left[|X_i|>\dfrac {x_i}{\theta_{i}},\,X_j >\dfrac {x_j}{\theta_{j}}\right] \PP\left[\Theta_{i}\in d\theta_{i}, \Theta_{j}\in d\theta_{j}\right],\\
	J_{22}(x_i,x_j) &:=& \int_{(b_i(x_i),\infty)} \int_{(b_j(x_j),\infty)} 
	\PP\left[|X_i|>\dfrac {x_i}{\theta_{i}},\,X_j >\dfrac {x_j}{\theta_{j}}\right] \PP\left[\Theta_{i}\in d\theta_{i}, \Theta_{j}\in d\theta_{j}\right].
	\eeao
	In a similar way, as for relations \eqref{eq.KLP.16}, \eqref{eq.KLP.17} and \eqref{eq.KLP.18}, we conclude that $J_{12}(x_i,\,x_j)$, $J_{21}(x_i,\,x_j)$ and  $J_{22}(x_i,\,x_j)$ are negligible with respect to $\PP[\Theta_i\,|X_i| >  x_i\,,\;\Theta_j\,X_j > x_j]$ as $x_i\wedge x_j \rightarrow \infty$. Hence, because $\PP[\Theta_i\,|X_i| >  x_i\,,\;\Theta_j\,X_j > x_j] \leq \PP[\Theta_j\,X_j > x_j]$, they are also negligible with respect to
	\beam \label{eq.KLP.22}
	\PP[\Theta_j\,X_j > x_j]\,,
	\eeam 
	as $x_i\wedge x_j \rightarrow \infty$  and 
	\beam \label{eq.KLP.23}
	J_{11}(x_i,\,x_j) \leq \PP\left[|X_i|>\dfrac {x_i}{b_{i}(x_i)},\,X_j >\dfrac {x_j}{b_{j}(x_j)}\right]= o\left(\PP[X_j > x_j]\right)\,,
	\eeam 
	as $x_i\wedge x_j \rightarrow \infty$, where in the last step we used the definition of the functions $b_i$, $b_j$ and the fact that the $X_i$, $X_j$ are $TAI$ and their distributions are in the class $\mathcal{D}$. Whence, from relations \eqref{eq.KLP.22}, \eqref{eq.KLP.23} in combination with \eqref{eq.KLP.21} we conclude
	\beam \label{eq.KLP.24} \notag
	\PP\left[ |\Theta_{i}\,X_i| > x_i\,,\;\Theta_j\,X_j > x_j\right]&\leq& o\left(\PP\left[X_j>x_j\right]\right)+o\left(\PP\left[\Theta_j\,X_j>x_j\right]\right)\\[2mm]
	&=&o\left(\PP\left[\Theta_j\,X_j>x_j\right]\right)\,,
	\eeam
	as $x_i\wedge x_j \rightarrow \infty$, where the last step follows from the fact that
	\beam \label{eq.KLP.25} 
	\PP\left[\Theta_j\,X_j>x_j\right]\asymp \PP\left[X_j>x_j\right]\,,
	\eeam
	as $x_j \rightarrow \infty$. 
	
	Indeed, on the one hand,  using Assumption \ref{ass.KLP.2} (see also Remark \ref{rem.KLP.2}) we obtain
	\beao
	\PP\left[\Theta_j\,X_j>x_j\right]&=&\left(\int_0^{b_j(x_j)} +\int_{b_j(x_j)}^{\infty}\right) \PP\left[X_j>\dfrac {x_j}{\theta_{j}}\right]\,\PP[\Theta_{j}\in d\theta_{j}] \\[2mm]
	&\leq&  \PP\left[X_j >\dfrac {x_j}{b_j(x_j)} \right]+\PP[\Theta_j>b_j(x_j)]\\
	&=&\PP\left[X_j >\dfrac {x_j}{b_j(x_j)} \right]+o\left( \PP[\Theta_j\,X_j >x_j]\right)\,,
	\eeao
	as $x_j \to \infty$. Hence, from property of class $\mathcal{D}$, and from the convergence $b_j(x_j) \to \infty$ as $x_j \to \infty$, we find
	\beam \label{eq.KP.3.a}
	\PP[\Theta_{j}\,X_j>x_j] \lesssim \PP\left[X_j > \dfrac{x_j}{b_j(x_j)}\right] \asymp \PP[X_j>x_j]\,,
	\eeam
	as $x_j \rightarrow \infty$. Therefore, by \eqref{eq.KP.3.a} it follows that
	\beam \label{eq.KP.3.b}
	\limsup_{x_j \to \infty} \dfrac{\PP[\Theta_{j}\,X_j>x_j]}{\PP[X_j>x_j]} < \infty \,.
	\eeam 
	On the other hand, for any arbitrarily chosen $\delta \in (0,\,1)$, we have the inequalities 
	\beam \label{eq.KP.3.c} \notag
	\PP[\Theta_{j}\,X_j>x_j] &\geq& \left(\int_{\delta}^1 + \int_1^{\infty} \right)\PP\left[X_j > \dfrac{x_j}{\theta_j}\right] \,\PP[\Theta_j \in d\theta_j] \\[2mm]
	&\geq& d\,\PP[X_j>x_j] \,\PP[\Theta_j \in (\delta, \,1]] + \PP[X_j>x_j] \,\PP[\Theta_j > 1]\\[2mm] \notag
	&\geq& (1\wedge d)\,\PP[X_j>x_j] \,\left(\PP[\Theta_j \in (\delta, \,1]] + \PP[\Theta_j > 1]\right) \longrightarrow d\,\PP[X_j>x_j]\,,
	\eeam 
	as $\delta \downarrow 0$, where $\PP\left[X_j > x_j/\theta_j \right] \geq \PP[X_j>x_j]$ in the second integral, that means when $\theta_j \in (1,\,\infty)$, while $d=d(\delta) \in (0,\,1]$ comes for any $\theta_j \in (\delta, \,1)$ from class $\mathcal{D}$. Thence by relation \eqref{eq.KP.3.c} we get
	\beao
	\liminf_{x_j \to \infty} \dfrac{\PP[\Theta_{j}\,X_j>x_j]}{\PP[X_j>x_j]} \geq\liminf_{x_j \to \infty} \dfrac{d\,\PP[X_j>x_j]}{\PP[X_j>x_j]} >0 \,,
	\eeao
	so we have
	\beam \label{eq.KP.3.d}
	\limsup_{x_j \to \infty} \dfrac{\PP[X_j>x_j] }{\PP[\Theta_j\,X_j >x_j]} < \infty\,.
	\eeam
	Whence, from \eqref{eq.KP.3.b} and \eqref{eq.KP.3.d} we get relation \eqref{eq.KLP.25}.
	
	Therefore, from \eqref{eq.KLP.24} and due to the arbitrariness of $i,\,j$ we obtain that $\Theta_1\,X_1,\,\ldots,\,\Theta_n\,X_n$ are $pTAI$.
	
	Now, we study the $pQAI$ case. We follow the same way, with the only difference that all the convergences hold for the same $x$. So, in relations \eqref{eq.KLP.22} and \eqref{eq.KLP.23} we find the terms
	\beao
	o\left(\PP\left[\Theta_i\,X_i>x\right]+\PP\left[\Theta_j\,X_j>x\right]\right),\qquad \qquad o\left(\PP\left[X_i>x\right]+\PP\left[X_j>x\right]\right), 
	\eeao
	respectively, as $x\to\infty$.
	~\halmos
	
	\smallskip
	
	Now we provide a dependence structure, which considers only dependencies between $\Theta_i$ and $X_i$ in each pair $(\Theta_i,\,X_i)$ for any $1\leq i \leq n$. This kind of dependence was introduced in \cite{asimit:jones:2008} via copulas, and studied extensively in risk theory, see \cite{yang:gao:li:2016} and \cite{yang:wang:leipus:siaulys:2013} for discrete-time models and \cite{asimit:badescu:2010} and \cite{li:tang:wu:2010} for continuous-time models.
	
	\begin{assumption} \label{ass.KLP.3} 
		For any pair $(X_i,\,\Theta_i)$, with $i\in\{1,\,\ldots,\,n\}$, there exists a function  ${h_i:\;[0,\,\infty) \longrightarrow (0,\,\infty)}$, such that  
		\beam \label{eq.KLP.31}
		\PP\left[ X_i> x \;\big|\; \Theta_{i}= \theta \right] \sim h_i(\theta)\, \PP[ X_{i}> x ]\,,
		\eeam 
		as $\xto$, uniformly for any $\theta \in s_{\Theta_i}$. 
	\end{assumption} 
	
The uniformity in relation \eqref{eq.KLP.31} is understood in the sense of
	\beao
	\lim_{\xto} \sup_{\theta \in s_{\Theta_i}} \left|\dfrac{\PP\left[ X_i> x \;\big|\; \Theta_{i}= \theta \right] }{h_i(\theta)\,\PP[X_i>x]} -1 \right|=0.
	\eeao

Before our first result, we need another one assumption that presents, in some sense, a variation of Assumption \ref{ass.KLP.3*}, which is not sensitive with respect to extreme negative values.

\begin{assumption} \label{ass.KLP.3.A} 
Let $\{X_1,\ldots,X_n\}$ be a collection of the  primary r.v.s, and let $\{\Theta_1,\ldots,\Theta_n\}$ be the collection of  the random weights. We suppose that the  Assumption  \ref{ass.KLP.3*} holds with some measurable functions $g_{ij}\;:\;[0,\,\infty)^2 \to (0,\,\infty)$ with $1\le i\ne j\le n$, and further there exist  constants $0<L \leq U < \infty$, such that  
\beam \label{eq.KLP.3.A} \notag
L\,g_{ij}(\theta_i,\,\theta_j)\,\PP[ X_{i}> x_i\,,\;X_{j}> x_j ] &\lesssim& \PP\left[ X_i> x_i\,,\;X_{j}> x_j\;\big|\; \Theta_{i}= \theta_i\,,\;\Theta_{j}= \theta_j \right]\\[2mm] 
&\lesssim&U\,g_{ij}(\theta_i,\,\theta_j)\,\PP[ X_{i}> x_i\,,\;X_{j}> x_j ]\,,
\eeam 
as $x_i\wedge x_j\to \infty$, for any $1\leq i \neq j \leq n$, uniformly with respect to $(\theta_i,\,\theta_j) \in s_{\Theta_i}\times s_{\Theta_j}$. 
\end{assumption} 

\bre \label{rem.KP.3.A}
We note that in class $\mathcal{D}$, for the distributions of primary r.v.s, under Assumption \ref{ass.KLP.3*}, and under some tail-balance conditions between the right and left distribution tail, Assumption \ref{ass.KLP.3.A} comes easily.
\ere

\bre\label{bb}
Let us suppose that the conditions of Lemma \ref{lem.KLP.1} hold, accepting instead of Assumption \ref{ass.KLP.3*}, the Assumption \ref{ass.KLP.3.A} with $n=2$. Following the same way $($with the necessary modifications after relation \eqref{eq.KLP.18*}$)$, we obtain  
\beam\notag 
&&L\,\int_{[0,\,\infty)}\int_{[0,\,\infty)} g_{12}(\theta_1,\,\theta_2)\,\PP\left[ X_{1}> \dfrac{x_1}{\theta_1}\,,\;X_{2}> \dfrac{x_2}{\theta_2} \right]\,\PP[\Theta_1 \in d\theta_1\,,\;\Theta_2 \in d\theta_2]\nonumber\\[2mm] 
&&\hspace{30mm}\lesssim \PP\left[ \Theta_1\,X_1> x_1\,,\;\Theta_2\,X_{2}> x_2\right]\nonumber\\[2mm] \notag
&&\lesssim U\,\int_{[0,\,\infty)}\int_{[0,\,\infty)} g_{12}(\theta_1,\,\theta_2)\,\PP\left[ X_{1}> \dfrac{x_1}{\theta_1}\,,\;X_{2}> \dfrac{x_2}{\theta_2} \right]\,\PP[\Theta_1 \in d\theta_1\,,\;\Theta_2 \in d\theta_2]\,,\nonumber
\eeam 
as $x_1\wedge x_2 \to \infty$.
\ere
	
	\bth \label{th.KLP.1}
	Let $X_1,\,\ldots,\,X_n$ be real-valued r.v.s, that are $pTAI$ (or $pQAI$), with distributions $F_1,\,\ldots,\,F_n \in \mathcal{D}$ respectively, and let $\Theta_1,\,\ldots,\,\Theta_n$ be non-negative, non-degenerate at zero r.v.s, where the  
	$\Theta_1,\,\ldots,\,\Theta_n$ and $X_1,\,\ldots,\,X_n$ satisfy Assumptions \ref{ass.KLP.2}, \ref{ass.KLP.3} and \ref{ass.KLP.3.A}. Then the products $\Theta_1 X_1,\,\ldots,\,\Theta_n X_n$ are also $pTAI$ (or $pQAI$ respectively) random variables.   
	\ethe
	
	\bre \label{rem.KLPS.3.5}
	Taking into account that Assumptions \ref{ass.KLP.3} and \ref{ass.KLP.3.A} include the independence as a special case, the last theorem presents in some sense partial extension of Theorem 2.2 in \cite{li:2013}. Here the word 'partial' means that it depends on the nature of the conditions for the variables $\Theta_1,\ldots, \Theta_n$ in both cases, either in Theorem \ref{th.KLP.1} or in Theorem 2.2 from \cite{li:2013}. For example, in the case of bounded from above $\Theta_1,\,\ldots,\,\Theta_n$, in both theorems, the conditions are satisfied. In addition, it should be noted that  in case the  $X_1,\,\ldots,\,X_n$ are non-negative r.v.s, all the results hold under Assumption \ref{ass.KLP.3*} instead of Assumption \ref{ass.KLP.3.A}.
	\ere
	
	\begin{proof}[\bf Proof of Theorem \ref{th.KLP.1}]
	We restrict ourselves only to the $pTAI$ case, as the proof for the $pQAI$ follows a similar way.
	
	For any $1\leq i \neq j \leq n$, let $\Theta_i^*,\,\Theta_j^*$ be new random variables, independent of any other sources of randomness, with the following joint distribution 
	\beao
	{\bf G}_{i j}^*(d\theta_i,\,d\theta_j):=g_{i j}(\theta_i,\,\theta_j)\,{\bf G}_{i j}(d\theta_i,\,d\theta_j)\,,
	\eeao
	where $g_{i j}$ is the function defined in Assumption \ref{ass.KLP.3*} and ${\bf G}_{i j}$ is the joint distribution of the random pair $(\Theta_i,\,\Theta_j)$. From relation \eqref{eq.KLP.13b} we have that ${\bf G}_{i j}^*$ represents a proper joint distribution. The corresponding marginal distributions of  $\Theta_i^*,\,\Theta_j^*$ are
	\beao
	G_{i}^{*}(d\theta_{i}):=h_{i}(\theta_{i})\,G_i(d\theta_{i}),\quad G_{j}^{*}(d\theta_{j}):=h_{j}(\theta_{j})\,G_j(d\theta_{j})\,.
	\eeao
	These marginals are also proper, see, e.g., \cite{chen:xu:cheng:2019} or \cite{li:tang:wu:2010}.  Furthermore,
	\beam \label{eq.KLP.27} \notag
	\PP[\Theta_i^*>b_i(x_i)]&=& \int_0^{\infty} \int_{b_i(x_i)}^\infty {\bf G}_{i j}^*(d\theta_i,\,d\theta_j)=\int_0^{\infty} \int_{b_i(x_i)}^\infty g_{i j}(\theta_i,\,\theta_j)\,{\bf G}_{i j}(d\theta_i,\,d\theta_j)\\[2mm]
	&\leq& K\int_0^{\infty} \int_{b_i(x_i)}^\infty {\bf G}_{i j}(d\theta_i,\,d\theta_j)=K\,\PP[\Theta_i> b_i(x_i)]\\[2mm] \notag
	&=&o\left( \PP\left[\Theta_i\,X_i>x_i\,,\;\Theta_j\,X_j>x_j\right] \right)\,,
	\eeam
	as $x_i \wedge x_j \to \infty$, because of Assumption \ref{ass.KLP.2}. Similarly, we obtain
	\beam \label{eq.KLP.28} 
	\PP[\Theta_j^*>b_j(x_j)]=o\left( \PP\left[\Theta_i\,X_i>x_i\,,\;\Theta_j\,X_j>x_j\right] \right)\,,
	\eeam
	as $x_i \wedge x_j \to \infty$. Now, by the assertion of Remark \ref{bb}, we have
	\beam \label{eq.KLP.29} 
	\PP\left[\Theta_i\,X_i>x_i\,,\;\Theta_j\,X_j>x_j\right] \asymp \PP\left[\Theta_i^*\,X_i>x_i\,,\;\Theta_j^*\,X_j>x_j\right]\,,
	\eeam
	as $x_i \wedge x_j \to \infty$. Hence, from relation \eqref{eq.KLP.29}, in combination with \eqref{eq.KLP.27} and \eqref{eq.KLP.28} we find
	\beao
	\PP\left[\Theta_i^*>b_i(x_i)\right] &=& o\left( \PP\left[\Theta_i^*\,X_i>x_i\,,\;\Theta_j^*\,X_j>x_j\right] \right)\,, \\
	\PP\left[\Theta_j^*>b_j(x_j)\right] &=&o\left( \PP\left[\Theta_i^*\,X_i>x_i\,,\;\Theta_j^*\,X_j>x_j\right] \right)\,,
	\eeao
	as $x_i \wedge x_j \to \infty$ implying that the primary  r.v.s $X_1,\ldots,X_n$ and random weights $\Theta_1^*,\ldots, \Theta_n^*$ satisfy Assumption \ref{ass.KLP.2}. In such a situation, by Lemma \ref{lem.KLP.2} we get 
    $$    \PP\left[\Theta_i^*\,|X_i|>x_i\,,\;\Theta_j^*\,X_j>x_j\right]=o[\PP(\Theta_{j}^{*}X_{j}>x_{j})],
    $$
    which together with \eqref{eq.KLP.29}  implies that
    $$    \PP\left[\Theta_i\,|X_i|>x_i\,,\;\Theta_j\,X_j>x_j\right]=o[\PP(\Theta_{j}^{*}X_{j}>x_{j})]\,,
    $$
	as $x_i \wedge x_j \to \infty$.  According to Lemma 3.1 in \cite{chen:xu:cheng:2019}  we have
	\beao	\PP(\Theta_{i}^{*}X_{i}>x_{i})\sim\PP(\Theta_{i}X_{i}>x_{i}),\qquad	\PP(\Theta_{j}^{*}X_{j}>x_{j})\sim\PP(\Theta_{j}X_{j}>x_{j})\,,
	\eeao
	as $x_i \rightarrow\infty$ and $x_j\rightarrow\infty$ respectively because  of Assumptions \ref{ass.KLP.2}, \ref{ass.KLP.3} and Remark \ref{rem.KLP.2}.
	Thus, by the last two relations, we conclude that the products $\Theta_1 X_1,\,\ldots,\,\Theta_n X_n$ are $pTAI$.
	\end{proof}
	
	\smallskip
	
	The next result can be found in \cite[Th.\ 2.5]{konstantinides:passalidis:2024} and gives the closedness with respect to the product of distribution classes $\mathcal{C}$, $\mathcal{D}\cap \mathcal{L}$ and $\mathcal{D}$, for the dependence structure of Assumption~\ref{ass.KLP.3}, under a sufficient but not necessary condition.
	
	\bpr \label{pr.KLP.2}
	Let $X$ be a real-valued r.v. with distribution $F$, $\Theta$ be non-negative, non-degenerate at zero r.v. with distribution $G$, and let $H$ be the distribution of the product $\Theta\,X$. Let hold Assumption \ref{ass.KLP.3} and for any constant $c>0$ it holds that $\bG(c\,x)=o[\bH(x)]$, as $\xto$. Then, the following statements hold
	\begin{enumerate}
		\item
		if $F \in \mathcal{D} \cap \mathcal{L}$, then $H \in \mathcal{D} \cap \mathcal{L}$;
		\item
		if $F \in \mathcal{C}$, then $H \in \mathcal{C}$;
		\item
		if $F \in \mathcal{D}$, then $H \in \mathcal{D}$.
	\end{enumerate}
	\epr  
	
	In the next result, the max-sum equivalence for randomly weighted sums is considered. A similar approach for the case of $pTAI$ can be found in \cite{li:2013}, and for the case of $pQAI$ can be found in \cite{chen:yuen:2009}. However, in both works, independence is assumed between the random weights and the primary random variables. Hence, as far the Assumption \ref{ass.KLP.3} contains the independence as a special case; our result contains these results as 'partial' cases.
	
	\bth \label{th.KLP.2}
	Let $X_1,\,\ldots,\,X_n$ be real-valued r.v.s,  with the corresponding distributions $F_1,\,\ldots,\,F_n \in \mathcal{D}\cap \mathcal{L}$ respectively, and let $\Theta_1,\,\ldots,\,\Theta_n$ be non-negative, non-degenerate at zero r.v.s. Under Assumptions \ref{ass.KLP.2}, \ref{ass.KLP.3} and \ref{ass.KLP.3.A}, the following statements hold:
	\begin{enumerate}
		\item
		If $X_1,\,\ldots,\,X_n$ are $pTAI$, then
		\beam \label{eq.KLP.33}
		\PP\left[ S_n^{\,\Theta}> x \right] \sim \PP\left[ \bigvee_{i=1}^n S_{i}^{\,\Theta}> x \right]\sim \PP\left[ \bigvee_{i=1}^n \Theta_{i}\,X_{i}> x \right]\sim \sum_{i=1}^n\PP[\Theta_{i}\,X_{i}> x]\,,
		\eeam
		as $\xto$. 
		\item
		If $X_1,\,\ldots,\,X_n$ are $pQAI$ and $F_1,\,\ldots,\,F_n \in \mathcal{C}$, then 
		\beam \label{eq.KLP.34}
		\PP\left[ S_n^{\Theta}> x \right] \sim \PP\left[ \bigvee_{i=1}^n S_{i}^{\Theta}> x \right]\sim \sum_{i=1}^n\PP[ \Theta_{i}\,X_{i}> x ]\,,
		\eeam
		as $\xto$. 
	\end{enumerate}
	\ethe
	
	\begin{proof}[\bf Proof]    
(1) From Proposition \ref{pr.KLP.2}(1) we observe that, for each $i=1,\dots,n$, the random variable $Z_i:=\Theta_i\,X_i$, follows the distribution $H_i \in \mathcal{D}\cap \mathcal{L}$, since relation \eqref{eq.KLP.13} is implied by Assumption \ref{ass.KLP.2} and Remark \ref{rem.KLP.2}. Hence, from 
		$\mathcal{D}\cap\mathcal{L} \subsetneq \mathcal{D}$, we can apply Theorem \ref{th.KLP.1} to find that the products $\Theta_1\,X_1,\,\ldots,\,\Theta_n\,X_n$ are also $pTAI$. Next, we use Proposition~\ref{pr.KLP.1}\,(1) to find relation \eqref{eq.KLP.33} after substitution of $Z_i$ by $\Theta_i\,X_i$. 
		        \smallskip
                
	(2)	By similar steps, using Proposition~\ref{pr.KLP.2}(2), the inclusion $\mathcal{C}\subsetneq \mathcal{D}$ and Proposition \ref{pr.KLP.1}(2), after substitution of $Z_i$ by $\Theta_i\,X_i$ we obtain
		\beao
		\PP\left[ S_n^{\Theta}> x \right] \sim \sum_{i=1}^n\PP[ \Theta_{i}\,X_{i}> x]\,,
		\eeao 
		as $\xto$. This asymptotic relation, in combination with the inequalities
		\beao
		\PP\left[ S_n^{\Theta}> x \right] \leq \PP\left[ \bigvee_{i=1}^n S_{i}^{\Theta}> x \right]\leq \PP\left[ \sum_{i=1}^n \Theta_{i}\,X_{i}^+> x \right],
		\eeao
		gives relation \eqref{eq.KLP.34}, which finishes the proof of the theorem. 
        \end{proof}
\smallskip	
	
	In the case of regular variation, we can get the direct expression under the additional condition of finite moments. In the next lemma, we generalize Breiman's lemma, in case of dependence described in Assumption \ref{ass.KLP.3}.
	
	\ble \label{lem.KLP.3.3}
	Let $X$ be a real-valued r.v. with distribution $F$, and $\Theta$ be a non-negative, non-degenerate at zero r.v. with distribution $G$, and let $H$ the distribution of the product $\Theta\,X$. If Assumption \ref{ass.KLP.3} holds, $\overline{G}(cx)=o\big(\PP(\Theta X>x)\big)$  for any constant $c>0$, $F\in \mathcal{R}_{-\alpha}$, for  $0<\alpha < \infty$, and $\E\left[ \Theta^{\rho}\,h(\Theta)\right] < \infty$, for some $\rho > \alpha$, then
	\beam \label{eq.KLP.a}
	\PP[\Theta\,X > x] \sim \E\left[ \Theta^{\alpha}\,h(\Theta)\right]\,\bF(x)\,,
	\eeam
	as $\xto$, which further implies  that $H\in \mathcal{R}_{-\alpha}$.
	\ele
	
	\pr~
	Let $\Theta^*$ be a random variable independent of any other source of randomness, with distribution $G^*(dt):=h(t)\,G(dt)$. Thus, we find that
	\beao
	\E\left[ \left(\Theta^*\right)^{\rho}\right]=\int_0^{\infty} \theta^{\rho}\,G^*(d\theta) = \int_0^{\infty} \theta^{\rho}\,h(\theta)\,G(d\theta) = \E\left[\Theta^{\rho}\,h(\Theta)\right] < \infty\,,
	\eeao
	for some $\rho>\alpha$, hence applying Breiman's lemma for the product $\Theta^*\,X$, we obtain
	\beam \label{eq.KLP.d}
	\PP\left[ \Theta^*\,X >x\right] \sim \E\left[\left(\Theta^*\right)^{\alpha}\right]\,\bF(x) =  \E\left[\Theta^{\alpha}\,h(\Theta)\right]\,\bF(x)\,,
	\eeam
	as $\xto$. Further, from Assumption \ref{ass.KLP.3},  relation $\overline{G}(cx)=o\big(\PP(\Theta X>x)\big)$, and Lemma 3.1 from \cite{chen:xu:cheng:2019} we find out that
	\beam \label{eq.KLP.e}
	\PP\left[ \Theta\,X >x\right] \sim \int_0^{\infty} h(\theta)\,\bF\left(\dfrac x\theta \right)\,G(d\theta) = \int_0^{\infty} \bF\left(\dfrac x\theta \right)\,G^*(d\theta) =\PP\left[\Theta^{*}\,X>x\right]\,,
	\eeam
	as $\xto$. So from relations \eqref{eq.KLP.d} and \eqref{eq.KLP.e} we have that \eqref{eq.KLP.a} is true. 
	Next, for any $b>0$ we get
	\beao
	\lim_{\xto} \dfrac{\bH(b\,x)}{\bH(x)} = \lim_{\xto} \dfrac{\PP[\Theta\,X > b\,x]}{\PP[\Theta\,X>x]} =\lim_{\xto}\dfrac{\E\left[ \Theta^{\alpha}\,h(\Theta)\right]\bF(b\,x)}{\E\left[ \Theta^{\alpha}\,h(\Theta)\right]\,\bF(x)}=b^{\alpha}\,,
	\eeao
	where the last step follows from $F \in \mathcal{R}_{-\alpha}$. Hence $H \in  \mathcal{R}_{-\alpha}$.
	~\halmos
	
	\smallskip
	
The next result follows directly from Theorem \ref{th.KLP.2} and Lemma \ref{lem.KLP.3.3}. In fact, with the additional condition of finite moment, we get the direct asymptotic expression in which we do not need to calculate the distribution of the product. 
	
	\bco \label{cor.KLP.3.1}
	Under the conditions of Theorem  \ref{th.KLP.2}, but restricting the distribution class of $X_1,\,\ldots,\,X_n$ into $F_i \in \mathcal{R}_{-\alpha_i}$ for $\alpha_i >0$, with $i=1,\,\ldots,\,n$ and when $\E\left[\Theta_i^{\rho_i}\,h_i(\Theta) \right]< \infty$, for some $\rho_i>\alpha_i$, with $i=1,\,\ldots,\,n$, we have the following.
	\begin{enumerate}
		\item
		If $X_1,\,\ldots,\,X_n$ are $pTAI$, then
		\beam \label{eq.KLP.f}
		\PP\left[ S_n^{\Theta}> x \right] \sim \PP\left[ \bigvee_{i=1}^n S_{i}^{\Theta}> x \right]\sim \PP\left[ \bigvee_{i=1}^n \Theta_{i}\,X_{i}> x \right]\sim \sum_{i=1}^n\E\left[ \Theta_{i}^{\alpha_i}\,h_i(\Theta_{i}) \right]\,\bF_i(x)\,,
		\eeam
		as $\xto$. 
		\item
		If $X_1,\,\ldots,\,X_n$ are $pQAI$, then 
		\beam \label{eq.KLP.g}
		\PP\left[ S_n^{\Theta}> x \right] \sim \PP\left[ \bigvee_{i=1}^n S_{i}^{\Theta}> x \right]\sim \sum_{i=1}^n\E\left[ \Theta_{i}^{\alpha_i}\,h_i(\Theta_{i}) \right]\,\bF_i(x)\,,
		\eeam
		as $\xto$. 
	\end{enumerate}  
	\eco
	
	\pr~
	Since relation \eqref{eq.KLP.13} is satisfied under Assumption  \ref{ass.KLP.2}, we can apply Lemma \ref{lem.KLP.3.3} for the products $\Theta_i\,X_i$, with $i=1,\,\ldots,\,n$, to find
	\beam \label{eq.KLP.h}
	\PP\left[ \Theta_{i}\,X_i> x \right] \sim \E\left[ \Theta_{i}^{\alpha_i}\,h_i(\Theta_{i}) \right]\,\bF_i(x)\,,
	\eeam
	as $\xto$. Next, from the inclusion $\mathcal{R}_{-\alpha_i}\subset \mathcal{D}\cap\mathcal{L}$, relations \eqref{eq.KLP.33} and \eqref{eq.KLP.h} we obtain \eqref{eq.KLP.f}. For the second part, we observe that $\mathcal{R}_{-\alpha_i}\subset \mathcal{C}$, and from relations \eqref{eq.KLP.34} and \eqref{eq.KLP.h} we obtain \eqref{eq.KLP.g}.
	~\halmos
	
	\smallskip
	
	The last theorem (and corollary) gives the finite time ruin probability for the discrete-time risk model with stochastic discount factors, in finite horizon. For example, we depict by $X_i$ the net losses in the $i$-th period and by $\Theta_i$ the corresponding stochastic discount factor. Then the $S_n^{\Theta}$ represents the discounted aggregate net-losses and, with $x$ as initial capital, the ruin probability during the $n$ first periods is equal to
	\beam \label{eq.KLP.35}
	\psi(x,\,n):=\PP\left[ \bigvee_{i=1}^n S_{i}^{\Theta}> x \right].
	\eeam
	Therefore, the next corollary follows from relation \eqref{eq.KLP.35}.
	
	\bco \label{cor.KLP.1}
			Under the conditions of Theorem \ref{th.KLP.2} either (1) or (2), it holds that
		\beao
		\psi(x,\,n)\mathop{\sim}\limits_{x\rightarrow\infty}  \sum_{i=1}^n\PP[\Theta_{i}\,X_{i}> x].
		\eeao	
	\\
		Under the conditions of Corollary \ref{cor.KLP.3.1}, either $(1)$ or $(2)$, we have that 
		\beao
		\psi(x,\,n) \mathop{\sim}\limits_{x\rightarrow\infty}\sum_{i=1}^n \E\left[ \Theta_{i}^{\alpha_i}\,h_i(\Theta_{i}) \right]\,\bF_i(x).
		\eeao
		 
		\eco

	\section{Generalized moments of randomly weighted sums} \label{sec.KP.3}

	In this section, we focus on asymptotic bounds of the quantity
	\beam \label{eq.KLP.37}
	\E\left[\varphi\left(S_n^{\Theta}\right)\,{\bf 1}_{\left\{S_n^{\Theta}>x\right\}}\right] = \E\left[\varphi\left(S_n^{\Theta}\right)\;\big|\;S_n^{\Theta}>x\right]\,\PP\left[ S_n^{\Theta}>x\right],
	\eeam
	where $\varphi$ is some function. The results present extensions of previous results and provide a continuation of Theorem \ref{th.KLP.2} for the case of $pQAI$ dependence structure for $F_1,\,\ldots,\,F_n \in \mathcal{D}$. Further, we find asymptotic bounds for the distribution tail of $S_n^{\Theta}$, namely $\PP\left[S_n^{\Theta}>x\right]$, which comes as consequence of the main theorem and the observation that if in relation \eqref{eq.KLP.37} we consider instead of $\varphi$ the unit function, then we obtain the tail $\PP\left[ S_n^{\Theta}>x\right]$. Although relation \eqref{eq.KLP.37} has no direct application in risk theory, it is relevant to risk management. We note that if replace $\varphi$ with the identical function, then relation \eqref{eq.KLP.37} becomes 
	\beao
	\E\left[S_n^{\Theta}\,{\bf 1}_{\left\{S_n^{\Theta}>x\right\}}\right]= \E\left[S_n^{\Theta}\;\big|\;S_n^{\Theta}>x\right]\,\PP\left[ S_n^{\Theta}>x\right]\,,
	\eeao
	where the expression $\E\left[S_n^{\Theta}\;\big|\;S_n^{\Theta}>x\right]$, is recognized as the Expected Shortfall, useful for the regulatory authorities to control the system risk. Similar results can be found in \cite{jaune:siaulys:2022} and \cite{liu:Yang:2021}.
	
	The generalized moments in randomly weighted sums attracted recently the attention in several papers, as  \cite{dirma:nakliuda:siaulys:2023}, \cite{dirma:paukstys:siaulys:2021}, \cite{jaune:ragulina:siaulys:2018} and \cite{leipus:paukstys:siaulys:2021}. In all of these works we observe a gradual improvement of the asymptotic estimates, however the random weights were supposed, are independent of the primary r.v.s $X_1,\,\ldots,\,X_n$, which follow some dependence structure as $pTAI$ or $pQAI$, while the random weights $\Theta_1,\,\ldots,\,\Theta_n$ were arbitrarily dependent.
	
	Here we propose a 'partial' generalization of \cite[Th.\ 3]{dirma:nakliuda:siaulys:2023}, considering dependence structure between the random weights and the primary r.v.s. In the next proposition we recall Theorem~2 of \cite{dirma:nakliuda:siaulys:2023}.
	
	\bpr \label{pr.KLP.3}
	Let $Z_1,\,\ldots,\,Z_n$ be real-valued r.v.s, with distributions $H_1,\,\ldots,\,H_n \in \mathcal{D}$ respectively, that are $pQAI$, and $\varphi\;:\;\bbr \rightarrow \bbr_+=[0,\,\infty)$ be a function, satisfying the following properties
	\begin{enumerate}
		\item
		The $\varphi$ is ultimately non-decreasing,
		\item
		The  $\varphi$ is ultimately sub-homogeneous, namely it holds that $\varphi(2\,x) \leq C\,\varphi(x)$ for some constant $C>0$ and some sufficiently large $x$,
		\item
		The  $\varphi$ is ultimately differentiable,
		\item
		For any integer $1\leq i \leq n$, the generalized moment $\E\left[\varphi\left(Z_i\right)\,{\bf 1}_{\left\{Z_i>x\right\}}\right]$ ultimately exists (it is finite).
	\end{enumerate}
	Then it holds that
	\beao
	\sum_{i=1}^n L_{H_i} \E\left[\varphi\left(Z_i\right) {\bf 1}_{\left\{Z_i>x\right\}}\right]&\lesssim& \E\left[\varphi\left(\sum_{i=1}^n Z_i\right) {\bf 1}_{\left\{\sum_{i=1}^n Z_i>x\right\}}\right]\\[2mm]
&\lesssim&\sum_{i=1}^n \dfrac 1{L_{H_i}}\,\E\left[\varphi\left(Z_i\right) {\bf 1}_{\left\{Z_i>x\right\}}\right]\,,
	\eeao
	as $\xto$.
	\epr 
	
	Next, we need a lemma before giving the main result of the section. Let us recall that $H_i(x)=\PP[\Theta_i,\,X_i \leq x]$, for any integer $1\leq i \leq n$.
	
	\ble \label{lem.KLP.3}
	Let $X$ be a real-valued r.v., with distribution $F\in \mathcal{D}$ and let $\Theta$ be a non-negative, non-degenerate at zero random variable, with distribution $G$. If Assumption \ref{ass.KLP.3} holds and 
	$\bG(c\,x)=o[\bH(x)]$, as $\xto$ for any $c>0$, then $H \in\mathcal{D}$ and $L_H\geq L_F$.
	\ele
	
	\pr~
	From Proposition \ref{pr.KLP.2}(3), we obtain $H \in \mathcal{D}$ and 
	\beao
&&\dfrac 1{L_{H}}=\lim_{v\uparrow 1}\limsup_{\xto}\dfrac{\PP[\Theta\,X >v\,x]}{\PP[\Theta\,X >x]} \\[2mm]
&&\leq \lim_{v\uparrow 1}\limsup_{\xto}\dfrac{\PP[\Theta\,X >v\,x\,,\;\Theta \leq b_1(x)]}{\PP[\Theta\,X >x]} + \lim_{v\uparrow 1}\limsup_{\xto}\dfrac{\PP[\Theta\,X >v\,x\,,\;\Theta > b_1(x)]}{\PP[\Theta\,X >x]}\\[2mm]
&&\leq \lim_{v\uparrow 1}\limsup_{\xto}\dfrac{\PP[\Theta\,X >v\,x\,,\;\Theta \leq b_1(x)]}{\PP[\Theta\,X >x\,,\;\Theta \leq b_1(x)]}  + o(1)= \lim_{v\uparrow 1}\limsup_{\xto}\dfrac{\int_0^{b_1(x)} h(t)\bF\left(\dfrac{v\,x}{t}\right)\,G(dt)}{\int_0^{b_1(x)} h(t)\bF\left(\dfrac{x}{t}\right)\,G(dt)}\\
&&\leq \lim_{v\uparrow 1}\limsup_{\xto}\sup_{0<\theta \leq b_1(x)}\dfrac{\bF\left(\dfrac{v\,x}{\theta}\right)}{\bF\left(\dfrac{x}{\theta}\right)}= \lim_{v\uparrow 1}\limsup_{\xto}\sup_{x/b_1(x) \leq z}\dfrac{\bF\left(v\,z\right)}{\bF\left(z\right)}\,,
	\eeao
	as $\xto$, where $b_1(x)$ represents a function satisfying Assumption \ref{ass.KLP.2} for the $\Theta$, and the $o(1)$ follows from:
	\beao
	\PP[\Theta\,X >b\,x\,,\;\Theta > b_1(x)] \leq \PP[\Theta > b_1(x)] = o(\PP[\Theta\,X >x])\,,
	\eeao
	as $\xto$, as comes from Assumption \ref{ass.KLP.2} through Remark \ref{rem.KLP.2}. Therefore, we conclude that
	\beao
	\dfrac 1{L_{H}} \leq \lim_{v \uparrow 1}\limsup_{\xto}\sup_{x/b_1(x) \leq z}\dfrac{\bF\left(v\,z\right)}{\bF\left(z\right)} = \dfrac 1{L_{F}}.~\halmos
	\eeao 
	
	\bth \label{th.KLP.3}
	Let $X_1,\,\ldots,\,X_n$ be real-valued r.v.s, with distributions $F_1,\,\ldots,\,F_n \in \mathcal{D}$ respectively, that are $pQAI$, and let $\Theta_1,\,\ldots,\,\Theta_n$ be non-negative, non-degenerate at zero r.v.s, with the $X_1,\,\ldots,\,X_n$ and $\Theta_1,\,\ldots,\,\Theta_n$ to satisfy Assumptions \ref{ass.KLP.2}, \ref{ass.KLP.3} and \ref{ass.KLP.3.A} for any integer $1\leq i \leq n$. We assume that there exists function $\varphi\;:\;\bbr \rightarrow \bbr_+=[0,\,\infty)$, from Proposition~\ref{pr.KLP.3}, satisfying the conditions (1), (2), (3) and additionally
	\begin{enumerate}
		\item[($4^*$)]
		For any integer $1\leq i \leq n$ the following generalized moment exists:
		\beao
		\E\left[\varphi\left(\Theta_i\,X_i\right)\,{\bf 1}_{\left\{\Theta_i\,X_i>x\right\}}\right]< \infty\,. 
		\eeao
	\end{enumerate}
	Then it holds that
	\beam \label{eq.KLP.40}
	\sum_{i=1}^n L_{F_i} \E\left[\varphi\left(\Theta_i X_i\right) {\bf 1}_{\left\{\Theta_i X_i>x\right\}}\right]\lesssim \E\left[\varphi\left(S_n^{\Theta}\right) {\bf 1}_{\left\{S_n^{\Theta}>x\right\}}\right]\lesssim \sum_{i=1}^n \dfrac 1{L_{F_i}}\E\left[\varphi\left(\Theta_i X_i\right) {\bf 1}_{\left\{\Theta_i X_i>x\right\}}\right],
	\eeam
	as $\xto$.
	\ethe
	
	\pr~
	From Proposition \ref{pr.KLP.2}\,(3) we find that 
	r.v.\ $Z_i=\Theta_i\,X_i$ follows a distribution $H_i \in \mathcal{D}$ for any integer $1\leq i \leq n$, as comes from Assumption \ref{ass.KLP.2} and Remark \ref{rem.KLP.2}. Further, from Theorem \ref{th.KLP.1} we obtain that the products $\Theta_1\,X_1,\,\ldots,\,\Theta_n\,X_n$ are $pQAI$. Hence, applying Proposition \ref{pr.KLP.3} for the $\Theta_i,X_i$ and Lemma \ref{lem.KLP.3}, relation  \eqref{eq.KLP.40} follows.
	~\halmos
	
	\bre \label{rem.KLP.4}
	If $F_i \in \mathcal{C}$, $i=1,\dots,n$, then the asymptotic bounds in relation \eqref{eq.KLP.40} become asymptotic equivalence. 
	
	Furthermore, the functions $\varphi(x)=1$ and $\varphi(x)=x$, for any $x\in \bbr$, satisfy the conditions (1), (2) and (3) in Proposition \ref{pr.KLP.3}. In the first case, with $\varphi(x)=1$, under the conditions of Theorem \ref{th.KLP.3}, taking into consideration equality \eqref{eq.KLP.37} we have
	\beam \label{eq.KLP.41}
	\sum_{i=1}^n L_{F_i} \PP\left[\Theta_i\, X_i>x\right]\lesssim \PP\left[S_n^{\Theta}>x\right]\lesssim \sum_{i=1}^n \dfrac 1{L_{F_i}}\PP\left[\Theta_i\,X_i>x\right]\,,
	\eeam
	as $\xto$. Here we can notice that if $F_i \in \mathcal{C}$ for any integer $1\leq i \leq n$, the result of relation \eqref{eq.KLP.41}, coincides with the asymptotic behavior of the tail of  $S_n^{\Theta}$, as appears in Theorem \ref{th.KLP.2}\,(2).
	
	In the second case, when $\varphi(x)=x$, under the assumptions of Theorem \ref{th.KLP.3} we obtain
	\beam \label{eq.KLP.42}
	\sum_{i=1}^n L_{F_i} \E\left[\Theta_i X_i\;\big|\;\Theta_i X_i>x\right]\,\dfrac{\PP\left[\Theta_i\, X_i>x\right]}{\PP\left[ S_n^{\Theta}>x\right]}&\lesssim& \E\left[S_n^{\Theta}\;\big|\;S_n^{\Theta}>x\right]\\[2mm]\notag
	&\lesssim& \sum_{i=1}^n \dfrac 1{L_{F_i}}\E\left[\Theta_i X_i\;\big|\;\Theta_i X_i>x\right]\dfrac{\PP\left[\Theta_i\, X_i>x\right]}{\PP\left[ S_n^{\Theta}>x\right]},
	\eeam
	as $\xto$. Relation \eqref{eq.KLP.42} provides asymptotic bounds or even asymptotic equivalence in $\mathcal{C}$, for the expected shortfall.
	\ere
	
	Now, we are interested in asymptotic bounds of the following expression $\E\left[\Theta_j\, X_j\, {\bf 1}_{\left\{S_n^{\Theta}>x\right\}}\right]$, for any integer $1\leq j \leq n$, which is connected with the Marginal Expected Shortfall:
	\beam \label{eq.KLP.43}
	\E\left[\Theta_i X_i\;\big|\;S_n^{\Theta}>x\right],
	\eeam
	see relevant works, studying the asymptotic behavior of relation \eqref{eq.KLP.43} in \cite{chen:liu:2022}, \cite{jaune:ragulina:siaulys:2018} and \cite{li:2022}. The following result can be found in \cite[Lem.\ 6, 7]{jaune:ragulina:siaulys:2018}.
	
	\bpr \label{pr.KLP.4}
	Let $Z_1,\,\ldots,\,Z_n$ be real-valued r.v.s, with distributions $H_1,\,\ldots,\,H_n \in \mathcal{D}$ respectively.  
	\begin{enumerate}
		\item
		If $Z_1,\,\ldots,\,Z_n$ are $pQAI$ and it holds that $\bH_i(x)=O[\bH_j(x)]$ as $\xto$, for any  integer $1\leq i \leq n$ and for some integer $1\leq j \leq n$, such that $\E[Z_j^+]< \infty$, then 
		\beao
		\E\left[Z_j\,{\bf 1}_{\left\{\sum_{i=1}^n Z_i>x\right\}}\right]\lesssim \E\left[Z_j\,{\bf 1}_{\left\{Z_j>x\right\}}\right]\,,
		\eeao
		as $\xto$.
		\item
		If $Z_1,\,\ldots,\,Z_n$ are $pTAI$ and it holds that $\bH_i(x)=O[\bH_j(x)]$ as $\xto$, for any  integer $1\leq i \leq n$ and for some integer $1\leq j \leq n$, such that $\E[|Z_j|]<\infty$, then 
		\beao
		L_{H_j}\E\left[Z_j\,{\bf 1}_{\left\{ Z_j>x\right\}}\right]\lesssim \E\left[Z_j\,{\bf 1}_{\left\{\sum_{i=1}^n Z_i>x\right\}}\right]\lesssim \E\left[Z_j\,{\bf 1}_{\left\{ Z_j>x\right\}}\right]\,,
		\eeao
		as $\xto$.
	\end{enumerate}
	\epr
	
	In the next result, we find dependence structure between random weights and primary random variables.
	\bco \label{cor.KLP.2}
	Let  $X_1,\,\ldots,\,X_n$ be real-valued r.v.s, with distributions $F_1,\,\ldots,\,F_n \in \mathcal{D}$ respectively, and let $\Theta_1,\,\ldots,\,\Theta_n$ be non-negative, non-degenerate at zero r.v.s. We assume that the $\Theta_1,\,\ldots,\,\Theta_n$ and $X_1,\,\ldots,\,X_n$ satisfy Assumptions \ref{ass.KLP.2}, \ref{ass.KLP.3} and \ref{ass.KLP.3.A}, with $\PP[\Theta_i\,X_i >x]=O(\PP[\Theta_j\,X_j>x])$, as $\xto$ for any integer $1\leq i \leq n$ and some integer  integer $1\leq j \leq n$, such that $\E[\Theta_j\,|X_j|]< \infty$.   
	\begin{enumerate}
		\item
		If $Z_1,\,\ldots,\,Z_n$ are $pQAI$, then 
		\beao
		\E\left[\Theta_j\,X_j\,{\bf 1}_{\left\{S_n^{\Theta}>x\right\}}\right]\lesssim \E\left[\Theta_j\,X_j\,{\bf 1}_{\left\{\Theta_j\,X_j>x\right\}}\right]\,,
		\eeao
		as $\xto$.
		\item
		If $Z_1,\,\ldots,\,Z_n$ are $pTAI$, then 
		\beao
		L_{F_j}\E\left[\Theta_j\,X_j\,{\bf 1}_{\left\{\Theta_j\,X_j>x\right\}}\right]\lesssim \E\left[\Theta_j\,X_j\,{\bf 1}_{\left\{S_n^{\Theta}>x\right\}}\right]\lesssim \E\left[\Theta_j\,X_j\,{\bf 1}_{\left\{ \Theta_j\,X_j>x\right\}}\right]\,,
		\eeao
		as $\xto$.
	\end{enumerate}
	\eco
	
	\begin{proof}[\bf Proof]
	(1)		From Assumptions \ref{ass.KLP.2} and \ref{ass.KLP.3}, taking into account Remark \ref{rem.KLP.2} and applying Proposition \ref{pr.KLP.2}(3), we obtain that $Z_i:=\Theta_i\,X_i$ follows distribution $H_i\in \mathcal{D}$, for any integer $1\leq i \leq n$. Further, from Theorem \ref{th.KLP.1} we find that the dependence structure $pQAI$ remains for $Z_1,\,\ldots,\,Z_n$. Hence, applying Proposition  \ref{pr.KLP.4}(1) we have the desired result.
		\smallskip\\
        (2)		Similarly, with the only difference that we apply Proposition \ref{pr.KLP.4}(2) instead of Proposition \ref{pr.KLP.4}(1) and then Lemma~\ref{lem.KLP.3}.
        \end{proof}

	\bre \label{rem.KLP.5}
	The formulas in Corollary \ref{cor.KLP.2} provide bounds for the marginal expected shortfall. Indeed, from the assumptions of part (2) we find
	\beao
	L_{F_j}\E\left[\Theta_j\,X_j\;\big|\;\Theta_j\,X_j>x\right]\,\dfrac{\PP\left[\Theta_i\, X_i>x\right]}{\PP\left[S_n^{\Theta}>x\right]}&\lesssim& \E\left[\Theta_j\,X_j\;\big|\;S_n^{\Theta}>x\right]\\[2mm] 
	&\lesssim& \E\left[\Theta_j\,X_j\;\big|\;\Theta_j\,X_j>x\right]\,\dfrac{\PP\left[\Theta_i\, X_i>x\right]}{\PP\left[S_n^{\Theta}>x\right]}\,,
	\eeao
	as $\xto$, while we have only upper bound in the case of $pQAI$. We notice that the marginal expected shortfall has bounds of the following form: the individual expected shortfall times the ratio of the tail of itself over the tail of the whole sum. The last can be considered as a contribution of the tail of the individual in the total tail. In the class $\mathcal{C}$, we find asymptotic equivalence relation for this risk measure.
	\ere

	\section{Randomly weighted and stopped sums} \label{sec.KP.5}

	Now we examine the asymptotic behavior of the randomized form of the sum in \eqref{eq.KLP.1}, namely
	\beam \label{eq.KLP.5.1}
	S_N^{\Theta}:=\sum_{i=1}^N \Theta_i \,X_i\,,
	\eeam 
	where $N$ represents a discrete random variable with values from $\bbn$. In particular, we can find the asymptotic relations for the non-weighted random sum $S_N:=\sum_{i=1}^N X_i$. For the studies of the asymptotic behavior of $S_N$ and its distribution class, see for example \cite{denisov:foss:korshunov:2010}, \cite{leipus:siaulys:2012} and \cite{sprindys:siaulys:2020}, among others. Besides the asymptotic behavior of $S_N$, one can investigate the randomly stopped maximum, as well the randomly stopped maximum of sums, whose definitions in the weighted case are as follows:
	\beao
		X_{(N)}^{\Theta}:=
		\begin{cases}
			0,& \text{if}\; N=0\\ 
			\bigvee_{i=1}^N \Theta_i \,X_i, &\text{if}\; N\geq 1
		\end{cases}\,, \qquad S_{(N)}^{\Theta} :=
		\begin{cases}
			0,& \text{if}\; N=0\\ 
			\bigvee_{i=1}^N S_i^{\Theta}, &\text{if}\; N\geq 1
		\end{cases}.
	\eeao 
	For a detailed survey on randomly stopped sums, maximums and maximums of sums, when the summands are heavy-tailed, in non-weighted case, we refer to \cite{karaseviciene:siaulys:2024} and \cite{leipus:siaulys:danilenko:karaseviciene:2024}. The applications of randomly stopped sums can be found in risk theory, since usually a portfolio has unknown number of claims over a concrete time interval. Further, the randomly stopped sums help to find asymptotic estimations of the ruin probability, see for example \cite[Sec.\ 4]{cheng:2015}. Inspired by the papers \cite{dindiene:leipus:2015}, \cite{olvera-cravioto:2012} and \cite{yang:leipus:siaulys:2015}, we consider a discrete-time risk model with the general weighted version of the randomly stopped sums. However, now we assume that the $\Theta_1,\,\Theta_2,\,\ldots$ are identically distributed as well as the $X_1,\,X_2,\,\ldots$ are.
	
	The conditions for the summands of randomly weighted stopped sums are found in the following assumption.
	
	\begin{assumption} \label{ass.KLP.5} 
		Let $\Theta,\,\Theta_1,\,\Theta_2,\,\ldots$ be non-negative, identically distributed r.v.s with common distribution $G$ and let $X,\,X_1,\,X_2,\,\ldots$ be identically distributed real-valued r.v.s with common distribution $F$. We assume that the discrete r.v. $N$ is non-negative and non-degenerate at zero and it has support bounded from above, namely there exists an integer $r_N$ such that $\PP[N>r_N]=0$. Finally, we suppose that $N$ is independent of $\Theta$ and $X$.
	\end{assumption} 
	
	In our assumptions we avoid to suppose independence among  $X_1,\,X_2,\,\ldots$, and in spite the restriction of the bounded support of the random variable $N$ and the identical distributions of the product $\Theta_i\,X_i$, we use $pTAI$ and  $pQAI$ dependence structure for the  $X_1,\,X_2,\,\ldots$ and arbitrary dependence among the $\Theta_1,\,\Theta_2,\,\ldots$. It is worth to notice that, because of the identical distributions among $X_1,\,X_2,\,\ldots$ and among $\Theta_1,\,\Theta_2,\,\ldots$, as described in Assumption~\ref{ass.KLP.5}, it holds the equality $b_1=b_2=b$ for the functions in Assumption~\ref{ass.KLP.2}, the equality $g_{ij}=g$ for any $1\leq i,\,j \leq 2$ for the functions in Assumption~\ref{ass.KLP.3.A} and  the equality $h_{i}=h$ for any $1\leq i \leq n$ for the functions in Assumption \ref{ass.KLP.3}.
	
	\bth \label{th.KLP.5.1}
	\begin{enumerate}		
		\item
		Under the conditions of Theorem \ref{th.KLP.2}(1) and Assumption \ref{ass.KLP.5} it holds that
		\beam \label{eq.KLP.5.3}
		\PP[S_N^{\Theta} > x] \sim \PP[S_{(N)}^{\Theta} > x] \sim \PP[X_{(N)}^{\Theta} > x] \sim \E[N]\,\PP[\Theta\,X>x]\,,
		\eeam 
		as $\xto$. Further, the distributions of $S_N^{\Theta}$, $S_{(N)}^{\Theta}$ and $X_{(N)}^{\Theta}$ belong to class $\mathcal{D}\cap \mathcal{L}$.
		
		\item
		Under the conditions of Theorem \ref{th.KLP.2}\,(2) and Assumption \ref{ass.KLP.5}, we obtain
		\beam \label{eq.KLP.5.4}
		\PP[S_N^{\Theta} > x] \sim \PP[S_{(N)}^{\Theta} > x] \sim \E[N]\,\PP[\Theta\,X>x]\,,
		\eeam 
		as $\xto$. Furthermore, the distributions of $S_N^{\Theta}$ and $S_{(N)}^{\Theta}$ belong to class $\mathcal{C}$.
	\end{enumerate}
	\ethe
	
	\pr~
	\begin{enumerate}
		
		\item
		Under Assumption \ref{ass.KLP.5}, via total probability formula we find
		\beam \label{eq.KLP.5.5} \notag
		\PP[S_N^{\Theta} > x] &=& \sum_{n=1}^{r_N} \PP[S_n^{\Theta} > x] \,\PP[N=n] \\[2mm]
		&\sim& \sum_{n=1}^{r_N} n\,\PP[\Theta\,X > x] \,\PP[N=n]= \E[N]\,\PP[\Theta\,X>x]\,,
		\eeam 
		as $\xto$, where in the second step, through the dominated convergence theorem, we used relation \eqref{eq.KLP.33}, keeping in mind the identically distributed $\Theta$ and $X$ by Assumption \ref{ass.KLP.5}. For the random maximum of sums, from the inequality
		\beam \label{eq.KLP.5.6}
		\PP[S_N^{\Theta} > x] \leq \PP[S_{(N)}^{\Theta} > x] \leq \PP[S_N^{\Theta+}> x]\,,
		\eeam 
		as $x>0$, where we denote $S_N^{\Theta+}:=\sum_{i=1}^N \Theta_i \, X_i^+$, we obtain through relation \eqref{eq.KLP.5.5} that
		\beam \label{eq.KLP.5.7}
		\PP[S_{(N)}^{\Theta} > x] \sim \E[N]\,\PP[\Theta\,X>x]\,,
		\eeam 
		as $\xto$. Furthermore, for the randomly stopped maximum, via relation \eqref{eq.KLP.33} we have
		\beam \label{eq.KLP.5.8} \notag
		\PP[X_{(N)}^{\Theta} > x] &=&\sum_{n=1}^{r_N} \PP\left[\bigvee_{i=1}^n \Theta_i\,X_i > x\right]\,\PP[N=n] \\[2mm]
		&\sim& \sum_{n=1}^{r_N} n\PP[\Theta\,X > x]\,\PP[N=n] =\E[N]\,\PP[\Theta\,X>x]\,,
		\eeam 
		as $\xto$. Relations \eqref{eq.KLP.5.5}, \eqref{eq.KLP.5.7} and \eqref{eq.KLP.5.8} imply \eqref{eq.KLP.5.3}. Further, by Proposition~\ref{pr.KLP.2}$(1)$, we get the membership $H \in \mathcal{D}\cap \mathcal{L}$ with $H(x):=\PP[\Theta\,X\leq x]$. Therefore, by \eqref{eq.KLP.5.3}, and taking into account the closure property of $\mathcal{D}\cap \mathcal{L}$ with respect to strong equivalence, see \cite[Pr.\ 3.10(i)]{leipus:siaulys:konstantinides:2023}, we find that the distributions of $S_N^{\Theta}$, $S_{(N)}^{\Theta}$ and $X_{(N)}^{\Theta}$ belong to class $\mathcal{D}\cap \mathcal{L}$. 
		
		\item
		In this case, relation \eqref{eq.KLP.5.5} follows in a similar way from relation \eqref{eq.KLP.34}, and by relations \eqref{eq.KLP.5.5} and \eqref{eq.KLP.5.6} we find \eqref{eq.KLP.5.7}. Further, from \eqref{eq.KLP.5.5} and \eqref{eq.KLP.5.7} we get \eqref{eq.KLP.5.4}, and by Proposition \ref{pr.KLP.2}$(2)$ we obtain $H \in \mathcal{C}$. Therefore, via \eqref{eq.KLP.5.4} and the closure property of class $\mathcal{C}$ with respect to strong equivalence (see \cite[Pr.\ 3.5(i)]{leipus:siaulys:konstantinides:2023}) we find that the distributions of $S_N^{\Theta}$ and $S_{(N)}^{\Theta}$  belong to class $\mathcal{C}$.~\halmos
	\end{enumerate}
	
	\bre \label{rem.KLPS.5.1}
	We observe that, as for randomly weighted finite sums, in the case of randomly stopped and weighted sum, the pairwise $QAI$ dependence structure fails to provide an asymptotic expression for the maximums. 
	\ere
	
	In the next result, we find more direct forms of relations \eqref{eq.KLP.5.3} and \eqref{eq.KLP.5.4}, reducing the distribution class to regular variation and adding a moment condition for r.v.\ $\Theta$, see Lemma~\ref{lem.KLP.3.3}. We note that the conditions $F_i\in \mathcal{R}_{-\alpha_i}$ and $\E[\Theta^{p_i}\,h_i(\Theta)] < \infty$, for some $p_i>\alpha_i$ of Corollary \ref{cor.KLP.3.1} are understood as $F \in \mathcal{R}_{-\alpha}$ and $\E[\Theta^p\,h(\Theta)] < \infty$ respectively, for some $p>\alpha$ due to Assumption \ref{ass.KLP.5}.
	
	\bco \label{cor.KLPS.5.1}
	\begin{enumerate}
		\item
		Under the conditions of Corollary \ref{cor.KLP.3.1}$(1)$ and Assumption \ref{ass.KLP.5}, we obtain
		\beam \label{eq.KLP.5.9}
		\PP[S_N^{\Theta} > x] \sim \PP[S_{(N)}^{\Theta} > x] \sim \PP[X_{(N)}^{\Theta} > x] \sim \E[N]\,\E[\Theta^{\alpha}\,h(\Theta)]\,\bF(x)\,,
		\eeam 
		as $\xto$, and further, the distributions of $S_N^{\Theta}$, $S_{(N)}^{\Theta}$ and $X_{(N)}^{\Theta}$ belong to class $\mathcal{R}_{-\alpha}$.
		
		\item
		Under the conditions of Corollary \ref{cor.KLP.3.1}$(2)$ and Assumption \ref{ass.KLP.5}, we obtain
		\beam \label{eq.KLP.5.10}
		\PP[S_N^{\Theta} > x] \sim \PP[S_{(N)}^{\Theta} > x] \sim \E[N]\,\E[\Theta^{\alpha}\,h(\Theta)]\,\bF(x)\,,
		\eeam 
		as $\xto$, and further, the distributions of $S_N^{\Theta}$ and $S_{(N)}^{\Theta}$ belong to class $\mathcal{R}_{-\alpha}$.
	\end{enumerate}
	\eco
	
	\pr~
	Relations \eqref{eq.KLP.5.9} and \eqref{eq.KLP.5.10} are implied by relations \eqref{eq.KLP.5.3} and \eqref{eq.KLP.5.4} respectively, in combination with relation \eqref{eq.KLP.a}. Next, since $F\in \mathcal{R}_{-\alpha}$, via \eqref{eq.KLP.5.9} and \eqref{eq.KLP.5.10}, the $S_N^{\Theta}$, $S_{(N)}^{\Theta}$ (and $X_{(N)}^{\Theta}$ in case $(1)$ only) follow distributions that belong to $\mathcal{R}_{-\alpha}$, because of the closure property of the class $\mathcal{R}_{-\alpha}$ with respect to strong equivalence, see \cite[Pr. 3.3(i)]{leipus:siaulys:konstantinides:2023}.
	~\halmos
	
	\smallskip
	
	Finally, we provide an application to risk theory, through the ruin probability in random time horizon. This kind of ruin probability is very popular and practical topic, see \cite{wang:su:yang:2024} and \cite{yang:leipus:siaulys:2015} among others. Following the lines of Section \ref{sec.KP.2}, the $X_i$ represent net losses, the $\Theta_i$ are discount factors and $x$ depicts the initial surplus, with the difference that now we examine a concrete time horizon, but the number of claims is unknown. Thus the ruin probability takes the form
	\beam \label{eq.KLP.5.11}
	\psi(x,\,N):=\PP[S_{(N)}^{\Theta}>x].
	\eeam
	From \eqref{eq.KLP.5.11}, we get the following result:
	
	\bco \label{cor.KLPS.5.2}
	\begin{enumerate}
		
		\item
		Under the conditions of Theorem \ref{th.KLP.5.1}, either of part $(1)$ or of part~$(2)$, we obtain
		\beao
		\psi(x,\,N) \sim \E[N]\,\PP[\Theta\,X>x]\,,
		\eeao
		as $\xto$.
		
		\item
		Under the conditions of Corollary \ref{cor.KLPS.5.1}, either of part $(1)$ or of part $(2)$, we obtain
		\beao
		\psi(x,\,N) \sim \E[N]\,\E[\Theta^{\alpha}\,h(\Theta)]\,\bF(x)\,,
		\eeao 
		as $\xto$.
	\end{enumerate}
	\eco

\end{document}